\documentclass[11pt]{jfm}
\usepackage[usenames,dvipsnames,svgnames,table]{xcolor}  
\usepackage{amsfonts}
\usepackage{fancyhdr}
\usepackage[parfill]{parskip} 
\usepackage{graphicx}
\usepackage{amssymb}
\usepackage{amsmath}
\usepackage{epstopdf}
\usepackage{natbib}
\usepackage{bm,footnpag}
\usepackage{layout}
\usepackage{caption}
\usepackage{subfig}
\usepackage{cleveref}
\usepackage{float}
\usepackage{xcolor,rotating,picture,lipsum}
\usepackage{wrapfig}
\usepackage{lscape}
\usepackage{rotating}
\usepackage{epstopdf}
\usepackage{multirow}
\usepackage{ifthen}
\usepackage{pgf}
\usepackage{pgffor}
\usepackage[inline]{asymptote} 
\usepackage{parskip}
\usepackage{csquotes}
\usepackage[nodisplayskipstretch]{setspace}
\usepackage{pgfplots}
\usepackage{tkz-euclide}
\usepackage{subfig}
\usepackage[toc]{appendix}
\usepackage[section]{placeins}
\usepackage[latin1]{inputenc}
\pgfplotsset{compat=1.16}

\ifCUPmtlplainloaded \else
  \checkfont{eurm10}
  \iffontfound
   \IfFileExists{upmath.sty}
      {\typeout{^^JFound AMS Euler Roman fonts on the system,
                   using the 'upmath' package.^^J}%
       \usepackage{upmath}}
      {\typeout{^^JFound AMS Euler Roman fonts on the system, but you
                   dont seem to have the}%
       \typeout{'upmath' package installed. JFM.cls can take advantage
                 of these fonts,^^Jif you use 'upmath' package.^^J}%
      }
  \else
  \fi
\fi

\ifCUPmtlplainloaded \else
  \checkfont{msam10}
  \iffontfound
    \IfFileExists{amssymb.sty}
      {\typeout{^^JFound AMS Symbol fonts on the system, using the
                'amssymb' package.^^J}%
       \usepackage{amssymb}%
         
       \let\ge=\geqslant  
      }{}
  \fi
\fi

\ifCUPmtlplainloaded \else
  \IfFileExists{amsbsy.sty}
    {\typeout{^^JFound the 'amsbsy' package on the system, using it.^^J}%
     \usepackage{amsbsy}}
    {}
\fi

\newsavebox{\astrutbox}
\sbox{\astrutbox}{\rule[-5pt]{0pt}{20pt}}

\allowdisplaybreaks[1]
\title[A class of exact solutions of the Navier-Stokes equations in $\mathbb{R}^3$ and $\mathbb{R}^4$]{A class of exact solutions of the Navier-Stokes equations in three and four dimensions}
\author[R. K. Michael Thambynayagam]%
{R. K. Michael Thambynayagam%
  \thanks{email: michael.thambynayagam@gmail.com}}
\affiliation{Retired scientist, Schlumberger,  Houston, Texas, USA}
 \numberwithin{equation}{section}
\begin{document}
\maketitle
\begin{abstract}
A few basic, intuitive, properties of the Navier-Stokes system of equations for incompressible fluid flows are discussed in this paper. We present a rephrased interpretation of the Navier-Stokes equation in a space having an arbitrary number of dimensions.  We then derive spatially periodic solutions for the velocity and pressure fields that span an unbounded domain in three and four dimensions, given a smooth solenoidal initial velocity vector field.  In these solutions all  velocity components depend non-trivially on all  coordinate directions.\\\\
\end{abstract}
\vspace*{-0.10 in}
\protect{\bfseries{2000 Math Subject Classification:}} 35Q30, 76D05
\\\\
\protect{\bfseries{ORCID:}} https://orcid.org/0000-0002-1778-7327\\\\
\protect{\bfseries{Journal of publication:}} European Journal of Mechanics - B/Fluids Volume 100, July--August 2023, Pages 12-20
\vspace*{-0.15 in}
\section{Introduction}
The Navier-Stokes equations are a set of partial differential equations that were developed by Claude-Louis \cite{Nav} and George Gabriel  \cite{Sto} to describe the motion of a Newtonian fluid, which can be either liquid or gas. These equations are essential for modeling fluid behavior in fluid dynamics, and they are used in many mathematical physics applications. Unfortunately, the equations are nonlinear and thus difficult to solve analytically. To gain a better understanding of the nonlinearity of these equations, one can refer to \cite{lin}. For those who need background knowledge on the subject,  \cite{bat}, \cite{pdnr}, \cite{con} and \cite{lad}'s monographs provide useful information.\\\\
Analytical or exact solutions are invaluable in understanding physical phenomena, as they provide a simple explicit expression of the behaviour in terms of well-defined functions. They are also essential for testing numerical solvers, ensuring accuracy and reliability of the solutions. Moreover, analytical solutions offer more insight into the problem than a numerical tabulation alone. \cite{bar}  has provided an insightful assessment of the significance and importance of analytical solutions to partial differential equations.\\\\
A general overview of analytical solutions to the Navier-Stokes equations is given by \cite{wan}, \cite{oka} and \cite{lan}. In a broard sense, analytical solutions to the Navier-Stokes equations are divided into two classes. The first class involves solutions where the nonlinearity is weakened or completely removed from the solution structure \cite{fox}. Examples of this simplified analysis include \cite{cou} and \cite{poi} flows. The second class is the Beltrami class of solutions, which involve nonzero nonlinear terms; \cite{bel}. Two-dimensional solutions to the Navier-Stokes equation for incompressible flows were developed by \cite{tay}, while unsteady analytical solutions involving all three Cartesian velocity components were presented by \cite{eth}. Recently, a triple-periodic fully three-dimensional analytical solution to the Navier-Stokes Equation for unsteady incompressible fluid flows was presented by \cite{ant}, extending Ethier and Steinman's work to a class of solutions for the velocity vector field in three dimensions.\\\\
Analytical solutions can be used to benchmark numerical Navier-Stokes solvers, as proposed in this paper. These solutions are extensible to higher dimensions, such as four or more, and can help identify and quantify discrepancies caused by numerical instabilities. Furthermore, they can provide insight into the underlying structure of physical systems, allowing for a deeper understanding of their behavior which can lead to further research and development in a variety of scientific disciplines. By exploring boundary conditions and turbulence, scientists can gain a better understanding of the physical world and develop solutions to complex problems. Additionally, higher dimensional analytical solutions are essential for developing robust numerical algorithms which can be used in a wide range of applications, such as those found in chemical engineering and biomedical engineering. Finally, these solutions can be used to stress-test Navier-Stokes solvers used to analyze the behavior of fluids in complex systems.\\\\
Antuono's solution structure is essentially the same as that of \cite{eth}, where the unsteady terms, $\frac{{\partial \bm{v}}}{{\partial t}}$, are equated to the viscous terms in the momentum equation. The velocity vector field $\bm{v}$ and the advective term $\left(\bm{v} \cdot \nabla\right)\bm{v}$ are represented through a stream function $\psi$ $\left(\bm{v}=\nabla \times \psi\right)$ and as the gradient of a scalar function respectively. These two equations are then decomposed into a linear equation for the stream function and a nonlinear equation for the pressure field. The linear equation for the stream function is solved analytically by separation of variables. Antuono describes the velocity vector field in a tours $\mathbb{T}^3=\left[0,L\right]^3$ in two parts, characterized by positive and negative helicity respectively. The velocity components were normalized by a reference velocity $v_r$ so that the average kinetic energy per unit of mass is $\left(\frac{v_r^2}{2}\right)$ at $t=0$. Finally, the pressure field was obtained in terms of a reference pressure and reference density from the velocity field.\\\\
In a previous article \cite{tha1}, a method was presented to derive analytical solutions by rephrasing the Navier-Stokes equation into three distinct terms associated with linear viscous forces, inertial forces, and external forces applied to the fluid. This paper applies the method to derive a time evolution analytical solution for the velocity vector field and pressure that span the entire unbounded domain. The mathematical expressions, particularly equations $\left(\ref{4.4}\right)$-$\left(\ref{4.30}\right)$ and $\left(\ref{5.5}\right)$-$\left(\ref{5.36}\right)$, are tedious and exhaustively boring to derive, but they are straightforward once understood. The method presented in this paper is an effective way of deriving analytical solutions for velocity vector fields and pressure in unbounded domains.
\vspace*{-0.20 in}
\section{The fundamental problem}
In a space having an arbitrary number of dimensions, $\mathbb{R}^n$, the fundamental equations of momentum and mass conservation for incompressible viscous fluid flow fields are represented as follows:
%
\begin{eqnarray}
\label{2.1}
\frac{{\partial v_i }}
{{\partial t}} + g_i= \kappa \Delta v_i -\frac{1}
{\rho }\frac{{\partial p}}
{{\partial x_i }} + f_i\quad {x \in \mathbb{R}^n},\quad t \geqslant 0\quad 
\end{eqnarray}
where $v_i = {v_i}\left( {x,t} \right)\,\, {x \in \mathbb{R}^n},\,\, i = 1,2,...,n,\,\,$ is the solenoidal velocity vector field and $p=p\left( {x,t} \right)$ the pressure of the fluid.
%
\begin{eqnarray}
\label{2.2}
g_i = {g_i}\left( {x,t} \right) = \sum_{j = 1}^n {v_j \frac{{\partial v_i }}
{{\partial x_j }}} 
\end{eqnarray}
is the nonlinear inertial force, $f_i=f_i\left( {x,t} \right)$ are the components of an externally applied force, $\rho$ is the constant density of the fluid,  $\kappa$ is the positive coefficient of kinematical viscosity and $\Delta  = \sum_{i = 1}^n {\frac{{\partial ^2 }}{{\partial x_i^{(2)} }}}$ is the Laplacian in the space variables.\\\\
The conservation of mass yields the solenoidal (incompressibility) condition
\begin{eqnarray}
\label{2.3}
\protect{\bf{div}}\,v = \sum\limits_{i = 1}^n {\frac{{\partial v_i }}
{{\partial x_i }}}=0 \qquad
x \in \mathbb{R}^n,\,\,t \geqslant 0
\end{eqnarray}
Conservation law also implies that the energy dissipation of a viscous fluid is bounded by the initial kinetic energy.  Therefore, if the externally applied force does no net work on the fluid, the kinetic energy of the solution should be finite. There exists a constant $E$ such that
%
\begin{eqnarray}
\label{2.4}
\int\limits_{\mathbb{R}^n } {\left| \bm{v} \right|} ^2 dx < E\qquad t \geqslant 0
\end{eqnarray}
The Navier-Stokes equation $\left(\ref{2.1}\right)$ must be solved forward in time $t\ge 0$, starting from an initial solenoidal velocity field 
\begin{eqnarray}
\label{2.5}
v_i\left( {x,0} \right) = v_i^{0} \qquad
x \in \mathbb{R}^n 
\end{eqnarray}
with the pressure evolving in time to maintain the incompressibility constraint $\left(\ref{2.3}\right)$. The superscript $0$ is used to denote the value of a function at time zero. 
%
%
\begin{eqnarray}
\label{2.6}
g_i\left( {x,0} \right)=g_i^0  = \sum_{j = 1}^n {v_j^0 \frac{{\partial v_i^0}}
{{\partial x_j }}} 
\end{eqnarray}
%
\section{Rephrasing the Navier-Stokes equation: $\mathbb{R}^n=\left\{ { - \infty  < x_i  < \infty ; \,\, i = 1,2,...,n} \right\}$}
In this section, in order to be perspicuous, we repeat the material from \cite{tha1} relating to the rephrasing of the Navier-Stokes equation and the resulting \protect{{\itshape{lemma}}. 
Assuming that the divergence and the linear operator can be commuted, the pressure field may be formally obtained by taking the divergence of  $\left(\ref{2.1}\right)$ as a solution of the Poisson equation, which is
%
\begin{eqnarray}
\label{3.1}
\Delta p = \rho \sum\limits_{i = 1}^n {\frac{{\partial \left( {f_i  - g_i } \right)}}
{{\partial x_i }}} 
\end{eqnarray}
Equation $\left(\ref{3.1}\right)$ is called the pressure Poisson equation. In a paper by \cite{gre}, the use of the pressure Poisson equation in solving the Navier-Stokes equation is discussed. It is important to note that while $\left(\ref{2.1}\right)$ and $\left(\ref{2.3}\right)$ lead to the pressure Poisson equation $\left(\ref{3.1}\right)$, the reverse does not always hold true; that is, $\left(\ref{2.1}\right)$ and $\left(\ref{3.1}\right)$ do not necessarily lead to $\left(\ref{2.3}\right)$. Therefore, when deriving solutions of the Navier-Stokes equations, it is essential to ensure that the velocity vector field remains solenoidal at all times, as defining pressure at the initial time independent of velocity would render the problem ill-posed.\\\\ The general solution of the Poisson equation $\left(\ref{3.1}\right)$ is given by 
%
%
\begin{eqnarray}
\label{3.2}
p\left( {x,t} \right) &=&- \frac{\rho }
{{2\pi }}\int\limits_{\mathbb{R}^2 } {{\rm P}\left( {y ,t} \right)\ln \left( {\frac{1}
{{\sqrt {\mathcal{P}_n\left( {x,y} \right)} }}} \right)\prod\limits_{j = 1}^2 {dy_j } },\,\, n=2, \nonumber\\
p\left( {x,t} \right) &=&- \frac{{\rho \Gamma \left( {\frac{n}
{2}} \right)}}
{{2\left( {n - 2} \right)\pi ^{\frac{n}
{2}} }}\int\limits_{\mathbb{R}^n } {\frac{{{\rm P}\left( {y,t} \right)}}
{{\left\{ {\mathcal{P}_n\left( {x,y} \right)} \right\}^{\frac{{n - 2}}
{2}} }}} \prod\limits_{j = 1}^n {dy_j },\,\, n\geqslant 3\,\,\, 
\end{eqnarray}
where $\Gamma \left( z \right) = \int_0^\infty  {e^{ - u} u^{z - 1} du}$\quad$\left[ {\Re z > 0} \right]$, is the Gamma function,
%
\begin{eqnarray}
\label{3.3}
{\rm P}\left( {x,t} \right) = \sum\limits_{j = 1}^n {\frac{{\partial \left( {f_j  - g_j } \right)}}
{{\partial x_j }}} 
\end{eqnarray}
and
%
%
\begin{eqnarray}
\label{3.4}
\mathcal{P}_n\left( {x,y} \right) = \sum\limits_{j = 1}^n {\left( {x_j  - y_j } \right)^2 } 
\end{eqnarray}
Differentiating $\left(\ref{3.2}\right)$ with respect to $x_i$ we obtain 
%
%
\begin{eqnarray}
\label{3.5}
\frac{{\partial p }}
{{\partial x_i }} = \frac{{\rho \Gamma \left( {\frac{n}
{2}} \right)}}
{{2\pi ^{\frac{n}
{2}} }}\int\limits_{\mathbb{R}^n } {\frac{{\left( {x_i  - y_i } \right){\rm P}\left( {y,t} \right)}}
{{\left\{ {\mathcal{P}_n\left( {x,y } \right)} \right\}^{\frac{n}
{2}} }}} \prod\limits_{j = 1}^n {dy_j },\,\, n\geqslant 2\qquad  
\end{eqnarray}
Substituting for $\frac{{\partial p}}{{\partial x_i }}$ in  $\left(\ref{2.1}\right)$ the following can be formulated:
%
\begin{eqnarray}
\label{3.6}
\frac{{\partial v_i }}
{{\partial t}} &=&\kappa \Delta v_i  -\frac{{\Gamma \left( {\frac{n}
{2}} \right)}}
{{2\pi ^{\frac{n}
{2}} }}\int\limits_{\mathbb{R}^n } {\frac{{\left( {x_i  - y_i } \right){\rm P}\left( {y,t} \right)}}
{{\left\{ {\mathcal{P}_n\left( {x,y} \right)} \right\}^{\frac{n}
{2}} }}} \prod\limits_{j = 1}^n {dy_j }+f_i - g_i,\quad n\geqslant 2,\,\, {x \in \mathbb{R}^n},\,\, t \geqslant 0\qquad  
\end{eqnarray}
The difficulty in solving the system of equations $\left(\ref{2.1}\right)-\left(\ref{2.3}\right)$ stems from the presence of the nonlinear term $g_i$.  We therefore rephrase the Navier-Stokes equation $\left(\ref{3.6}\right)$ as:
%
\begin{eqnarray}
\label{3.7}
\frac{{\partial v_i }}
{{\partial t}} &=&\kappa \Delta v_i  - \mathcal{U}_i+\mathcal{F}_i\quad {x \in \mathbb{R}^n},\qquad t \geqslant 0\quad
\end{eqnarray}
%
where 
%
\begin{eqnarray}
\label{3.8}
\mathcal{U}_i  = g_i  - \frac{{\Gamma \left( {\frac{n}
{2}} \right)}}
{{2\pi ^{\frac{n}
{2}} }}\int\limits_{\mathbb{R}^n } {\frac{{\left( {x_i  - y_i } \right)\sum\limits_{k= 1}^n {\frac{{\partial g_k \left( {y,t} \right)}}
{{\partial y_k }}} }}
{{\left\{ {\mathcal{P}_n\left( {x,y} \right)} \right\}^{\frac{n}
{2}} }}} \prod\limits_{j = 1}^n {dy_j },\quad {x \in \mathbb{R}^n},\,\, t \geqslant 0 \qquad
\end{eqnarray}
and
\begin{eqnarray}
\label{3.9}
\mathcal{F}_i= f_i  - \frac{{\Gamma \left( {\frac{n}
{2}} \right)}}
{{2\pi ^{\frac{n}
{2}} }}\int\limits_{\mathbb{R}^n } {\frac{{\left( {x_i  - w_i } \right)\sum\limits_{k = 1}^n {\frac{{\partial f_k \left( {w ,t} \right)}}
{{\partial w_k }}} }}
{{\left\{ {\mathcal{P}_n\left( {x,w} \right)} \right\}^{\frac{n}
{2}} }}} \prod\limits_{j = 1}^n {dw_j },\quad {x \in \mathbb{R}^n},\,\, t \geqslant 0 \qquad
\end{eqnarray}
$\mathcal{U}_i=\mathcal{U}_i\left( x,t \right)$ and $\mathcal{F}_i=\mathcal{F}_i\left( {x ,t} \right)$. The three terms on the right hand side of $\left(\ref{3.7}\right)$, $\kappa \Delta v_i $,  $\mathcal{U}_i$ and $\mathcal{F}_i$ are associated, respectively, with the linear viscous force,  the nonlinear inertial force and the externally applied force acting on the fluid.\\\\
A posteriori we state the following \protect{{\itshape{lemma}}: If $\mathcal{U}_i \equiv 0$, that is
%
\begin{eqnarray}
\label{3.10}
g_i = \frac{{\Gamma \left( {\frac{n}
{2}} \right)}}
{{2\pi ^{\frac{n}
{2}} }}\int\limits_{\mathbb{R}^n } {\frac{{\left( {x_i  - y_i } \right)\sum\limits_{k = 1}^n {\frac{{\partial g_k \left( {y,t} \right) }}
{{\partial y_k }}} }}
{{\left\{ {\mathcal{P}_n\left( {x,y} \right)} \right\}^{\frac{n}
{2}} }}} \prod\limits_{j = 1}^n {dy_j },
\end{eqnarray}
then the solution of the Navier-Stokes equation $\left(\ref{2.1}\right)$ is determined as a solution of the non-homogeneous diffusion equation \cite{car, tha2}:
%
\begin{eqnarray}
\label{3.11}
v_i\left( {x,t} \right) &=&\frac{1}
{{\left( {2\sqrt {\pi \kappa t} } \right)^n }}\int\limits_{\mathbb{R}^n } {v_i^{0} \left( {y} \right)e^{ - \sum\limits_{k = 1}^n {\frac{{\left( {x_k - y_k } \right)^2 }}
{{4\kappa t}}} } \prod\limits_{j = 1}^n {dy_j } }  + \nonumber\\
&+&\frac{1}{\left( {2\sqrt {\pi \kappa } } \right)^n }\int\limits_{\mathbb{R}^n } {\int\limits_0^t {\frac{{\mathcal{F}_i \left( {y ,\tau } \right)e^{ - \sum\limits_{k = 1}^n {\frac{{\left( {x_k  - y_k } \right)^2 }}
{{4\kappa \left( {t - \tau } \right)}}} } }}
{{\left( {t - \tau } \right)^{\frac{n}{2}}}}d\tau\prod\limits_{j = 1}^n {dy_j } } }
\end{eqnarray}
belong to the class of Beltrami flows. If the externally applied force on the fluid $f_i$ is set to zero, then, the second term on the right hand side of equation $\left(\ref{3.11}\right)$ vanishes and the solution of the Navier-Stokes equation $\left(\ref{2.1}\right)$ reduces to that of the Cauchy diffusion equation:
%
\begin{eqnarray}
\label{3.12}
v_i\left( {x,t} \right) &=&\frac{1}
{{\left( {2\sqrt {\pi \kappa t} } \right)^n }}\int\limits_{\mathbb{R}^n } {v_i^{0} \left( {y} \right)e^{ - \sum\limits_{k = 1}^n {\frac{{\left( {x_k - y_k } \right)^2 }}
{{4\kappa t}}} } \prod\limits_{j = 1}^n {dy_j } }  
\end{eqnarray}
We may express (3.5) as follows if $\left(\ref{3.10}\right)$ is satisfied for an $n\text{-tuple-periodic}$ form of the solenoidal velocity vector field: 
%
\begin{eqnarray}
\label{3.13}
\frac{{\partial p }}
{{\partial x_i }} =- \rho g_i=-\rho g_i^0e^{ - 2n\alpha ^2 \kappa t} 
\end{eqnarray}
Where $\alpha$ is the wavelength. Pressure is obtained from straightforward integration of $\left(\ref{3.13}\right)$.\\\\
Appendix A provides two known solutions, the two-dimensional double-periodic solution of \cite{tay}  and the unsteady unit periodic three-dimensional solution of \cite{arn} derived by \cite{tha1}, that satisfy the integral equation $\left(\ref{3.10}\right)$.
\vspace*{-0.50 in}
\section{A solution in  $\mathbb{R}^3=\left\{ { - \infty  < x_i  < \infty ; \,\, i = 1,2,3} \right\}$}
At time zero, we consider the solenoidal velocity vector field $v_i^{0}\left( x \right)$ of a triple-periodic form 
\begin{eqnarray}
\label{4.1}
\frac{v_1^0}{v_r}& = &\sin { \left(\alpha x_1+\xi_1 \right) } \cos {\left(\alpha x_2+\xi_2\right) }\sin { \left(\alpha x_3+\xi_3\right) } -\nonumber\\
&-&\sin  {\left( \alpha x_1+\xi_2\right) } \cos  { \left(\alpha x_3+\xi_1 \right) } \sin  {\left( \alpha x_2+\xi_3\right) }
\end{eqnarray}
%
\begin{eqnarray}
\label{4.2}
\frac{v_2^0}{v_r}& = &\sin { \left(\alpha x_2+\xi_1 \right) } \cos {\left(\alpha x_3+\xi_2\right) }\sin { \left(\alpha x_1+\xi_3\right) }-\nonumber\\
&-&\sin  {\left( \alpha x_2+\xi_2\right) } \cos  { \left(\alpha x_1+\xi_1 \right) } \sin  {\left( \alpha x_3+\xi_3\right) } 
\end{eqnarray}
%
\begin{eqnarray}
\label{4.3}
\frac{v_3^0}{v_r} &=& \sin { \left(\alpha x_3+\xi_1 \right) } \cos {\left(\alpha x_1+\xi_2\right) }\sin { \left(\alpha x_2+\xi_3\right) } -\nonumber\\
&-&\sin  {\left( \alpha x_3+\xi_2\right) } \cos  { \left(\alpha x_2+\xi_1 \right) } \sin  {\left( \alpha x_1+\xi_3\right) } 
\end{eqnarray}
where $v_r$ is the reference velocity, $\xi_1$, $\xi_2$, and $\xi_3$ are phase angles. The externally applied force on the fluid $f_i$ is set to zero. Equation $\left(\ref{3.10}\right)$ may be written as follows at time zero in $\mathbb{R}^3$:
%
\begin{eqnarray}
\label{4.4}
g_i^{0}= \frac{1}{{4\pi }}\int\limits_{\mathbb{R}^3 } {\frac{{\left( {x_i  - y_i } \right)\mathcal{G}^0\left(y_1, y_2,y_3\right) }}
{{\left\{ {\mathcal{P}_3\left( {x,y} \right)} \right\}^{\frac{3}
{2}} }}} \prod\limits_{j = 1}^3 {dy_j }\qquad i = 1,2,3
\end{eqnarray}
where $\mathcal{G}^0\left(y_1, y_2,y_3\right)=\sum_{k = 1}^3 \frac{\partial g_k^0 \left( {y,t} \right)}{\partial y_k }$. Changing the variable of integration, $u_i=x_i  - y_i$, gives
%
\begin{eqnarray}
\label{4.5}
g_i^{0}= \frac{1}{{4\pi }}
\int\limits_{\mathbb{R}^3 } {\frac{{u_i\,\mathcal{G}^0\left(x_1-u_1, x_2-u_2,x_3-u_3\right) }}
{{\left[u_1^2+u_2^2+u_3^2 \right]^{\frac{3}
{2}} }}}  \prod\limits_{j = 1}^3 {du_j }\quad i = 1,2,3
\end{eqnarray}
It becomes apparent on closer examination of $\left(\ref{4.5}\right)$ that any term in $g_{i}^0$, derived from $\left(\ref{2.6}\right)$, that is not a function of $x_i$ cannot be recovered by performing the integrals on the righthand side of $\left(\ref{4.5}\right)$; as a result of the integral identity $\int_{ - \infty }^\infty  {\frac{u}{u^2+\beta}} du = 0$, $\beta$ a constant, these terms vanish entirely. As a consequence, it is a prerequisite for the integral equation  $\left(\ref{4.5}\right)$ to hold that all of the terms that are not functions of $x_i$ on the righthand side of $\left(\ref{2.6}\right)$ must sum to zero. We therefore divide $g_i^0$ into two parts as follows:
%
\begin{eqnarray}
\label{4.6}
g_i^0 = \mathcal{V}_i^0 + \mathcal{W}_i^0\qquad i = 1,2,3
\end{eqnarray}
$\mathcal{V}_i^0$ is the sum of all terms that are not functions of $x_i$, and $\mathcal{W}_i^0$ is the sum of all remaining terms. For the integral equation  $\left(\ref{4.5}\right)$ to hold all terms in $\mathcal{V}_i^0$ must sum to zero. By substituting $v_1^0$, $v_2^0$ and $v_3^0$ into $\left(\ref{2.6}\right)$ and separating the terms into  $\mathcal{V}_i^0$ and $\mathcal{W}_i^0$ we get the following expressions:
%
\begin{eqnarray}
\label{4.7}
g_{1}^0\left( {x_1,x_2,x_3; \xi_1, \xi_2, \xi_3} \right)&=&\mathcal{V}_1^0\left( {x_2,x_3; \xi_1, \xi_2, \xi_3} \right) + \mathcal{W}_1^0\left( {x_1,x_2,x_3; \xi_1, \xi_2, \xi_3} \right)
\end{eqnarray}
%
\begin{eqnarray}
\label{4.8}
g_{2}^0\left( {x_2,x_3,x_1; \xi_1, \xi_2, \xi_3} \right)&=&\mathcal{V}_2^0\left( {x_3,x_1; \xi_1, \xi_2, \xi_3} \right) +  \mathcal{W}_2\left( {x_2,x_3,x_1; \xi_1, \xi_2, \xi_3} \right)
\end{eqnarray}
%
\begin{eqnarray}
\label{4.9}
g_{3}^0\left( {x_3,x_1,x_2; \xi_1, \xi_2, \xi_3} \right)&=&\mathcal{V}_3^0\left( {x_1,x_2; \xi_1, \xi_2, \xi_3} \right) +  \mathcal{W}_3\left( {x_3,x_1,x_2; \xi_1, \xi_2, \xi_3} \right)
\end{eqnarray}
%
where
%
\begin{eqnarray}
\label{4.10}
\frac{\mathcal{V}_1^0\left( {x_2,x_3; \xi_1, \xi_2, \xi_3} \right)}{v_r^2}&=&-\frac{\alpha}{4}\cos\left(\xi_1-\xi_3\right)\sin\left(\xi_2-\xi_3\right)\cos\left(2\alpha x_2+\xi_2+\xi_1\right)-\nonumber\\
&-&\frac{\alpha}{4}\cos\left(\xi_1-\xi_3\right)\cos\left(\xi_1-\xi_2\right)\sin\left(2\alpha x_3+\xi_2+\xi_3\right)-\nonumber\\
&-&\frac{\alpha}{4}\cos\left(\xi_3-\xi_2\right)\cos\left(\xi_1-\xi_2\right)\sin\left(2\alpha x_2+\xi_3+\xi_1\right)-\nonumber\\
&-&\frac{\alpha}{4}\cos\left(\xi_3-\xi_2\right)\sin\left(\xi_1-\xi_3\right)\cos\left(2\alpha x_3+\xi_1+\xi_2\right)+\nonumber\\
&+&\frac{\alpha}{4}\sin\left(\xi_2-\xi_1\right)\sin\left(\xi_3-\xi_1\right)\sin\left(2\alpha x_2+\xi_2+\xi_3\right)+\nonumber\\
&+&\frac{\alpha}{4}\sin\left(\xi_2-\xi_1\right)\sin\left(\xi_2-\xi_3\right)\sin\left(2\alpha x_3+\xi_3+\xi_1\right)
\end{eqnarray}
%
\begin{eqnarray}
\label{4.11}
\mathcal{V}_2^0\left( {x_3,x_1; \xi_1, \xi_2, \xi_3} \right)&=&v_r^2\mathcal{V}_1^0\left( {x_3,x_1; \xi_1, \xi_2, \xi_3} \right)
\end{eqnarray}
%
%
\begin{eqnarray}
\label{4.12}
\mathcal{V}_3^0\left( {x_1,x_2; \xi_1, \xi_2, \xi_3} \right)&=&v_r^2\mathcal{V}_2^0\left( {x_1,x_2; \xi_1, \xi_2, \xi_3} \right)
\end{eqnarray}
%
\begin{eqnarray}
\label{4.13}
\frac{\mathcal{W}_1^0\left( {x_1,x_2,x_3; \xi_1, \xi_2, \xi_3} \right)}{{v_r^2}}&=&\frac{\alpha }{4}\sin\left(2\alpha x_1+2\xi_1 \right)-\frac{\alpha }{4}\sin\left(2\alpha x_1+2\xi_1 \right)\cos \left(2\alpha x_3+2\xi_3 \right)+\nonumber\\
&+&\frac{\alpha }{4} \sin\left(2 \alpha x_1+2\xi_2\right)-\frac{\alpha }{4} \sin\left( 2\alpha x_1+2\xi_2\right)\cos\left(2\alpha x_2+2\xi_3 \right)-\nonumber\\
&-&\frac{\alpha}{2}\sin\left(\xi_3-\xi_2\right)\sin\left(\xi_3-\xi_1\right)\sin\left(2\alpha x_1+\xi_1+\xi_2\right)-\nonumber\\
&-&\frac{\alpha}{4}\sin\left(\xi_3-\xi_1\right)\sin\left(2\alpha x_1+\xi_1+\xi_2\right)\sin\left(2\alpha x_2+\xi_2+\xi_3\right)-\nonumber\\
&-&\frac{\alpha}{4}\sin\left(\xi_3-\xi_2\right)\sin\left(2\alpha x_1+\xi_1+\xi_2\right)\sin\left(2\alpha x_3+\xi_3+\xi_1\right)+\nonumber\\
&+&\frac{\alpha}{4}\sin\left(\xi_2-\xi_3\right)\cos\left(2\alpha x_1+\xi_3+\xi_1\right)\cos\left(2\alpha x_2+\xi_1+\xi_2\right)+\nonumber\\
&+&\frac{\alpha}{4}\cos\left(\xi_1-\xi_2\right)\cos\left(2\alpha x_1+\xi_3+\xi_1\right)\sin\left(2\alpha x_3+\xi_2+\xi_3\right)-\nonumber\\
&+&\frac{\alpha}{4}\cos\left(\xi_1-\xi_2\right)\cos\left(2\alpha x_1+\xi_2+\xi_3\right)\sin\left(2\alpha x_2+\xi_3+\xi_1\right)+\nonumber\\
&+&\frac{\alpha}{4}\sin\left(\xi_1-\xi_3\right)\cos\left(2\alpha x_1+\xi_2+\xi_3\right)\cos\left(2\alpha x_3+\xi_1+\xi_2\right)
\end{eqnarray}
%
\begin{eqnarray}
\label{4.14}
\mathcal{W}_2^0\left( {x_2,x_3,x_1; \xi_1, \xi_2, \xi_3} \right)&=&v_r^2\mathcal{W}_1^0\left( {x_2,x_3,x_1; \xi_1, \xi_2, \xi_3} \right)
\end{eqnarray}
%
\begin{eqnarray}
\label{4.15}
\mathcal{W}_3^0\left( x_3,{x_1,x_2; \xi_1, \xi_2, \xi_3} \right)&=&v_r^2\mathcal{W}_2^0\left( x_3,{x_1,x_2; \xi_1, \xi_2, \xi_3} \right)
\end{eqnarray}
If $\left(\ref{4.5}\right)$ is to hold, it is imperative that the following prerequisites are met:
\begin{eqnarray}
\label{4.16}
\mathcal{V}_1^0\left( {x_2,x_3; \xi_1, \xi_2, \xi_3} \right)=0\qquad \forall -\infty<x_2<\infty\,\,\textnormal{and}\,\,-\infty<x_3<\infty
\end{eqnarray}
%
\begin{eqnarray}
\label{4.17}
\mathcal{V}_2^0\left( {x_3,x_1; \xi_1, \xi_2, \xi_3} \right)=0\qquad \forall -\infty<x_3<\infty\,\,\textnormal{and}\,\,-\infty<x_1<\infty
\end{eqnarray}
%
\begin{eqnarray}
\label{4.18}
\mathcal{V}_3^0\left( {x_1,x_2; \xi_1, \xi_2, \xi_3} \right)=0\qquad \forall -\infty<x_1<\infty\,\,\textnormal{and}\,\,-\infty<x_2<\infty
\end{eqnarray}
It is therefore necessary to choose the phase angles $\xi_1$, $\xi_2$ and $\xi_3$ such that they satisfy the prerequisites $\left(\ref{4.16}\right)$, $\left(\ref{4.17}\right)$ and $\left(\ref{4.18}\right)$. In this particular case, as can be seen by substituting into $\left(\ref{4.10}\right)$, $\left(\ref{4.11}\right)$ and $\left(\ref{4.12}\right)$, the phase angles $\xi_1=-\frac{\pi}{3}$, $\xi_2=\frac{\pi}{3}$, and $\xi_3=\frac{\pi}{2}$ satisfy $\left(\ref{4.16}\right)$, $\left(\ref{4.17}\right)$ and $\left(\ref{4.18}\right)$.
%
\begin{eqnarray}
\label{4.19}
\frac{\mathcal{V}_1^0\left( {x_2,x_3;-\frac{\pi}{3}, \frac{\pi}{3}, \frac{\pi}{2}} \right)}{v_r^2}
&=&-\frac{\alpha\sqrt{3}}{16}\cos\left(2\alpha x_2\right)+\frac{3\alpha}{32}\sin\left(2\alpha x_3\right)-\frac{\alpha\sqrt{3}}{32}\cos\left(2\alpha x_3\right)+\nonumber\\
&+&\frac{3\alpha}{32}\sin\left(2\alpha x_2\right)+\frac{\alpha\sqrt{3}}{32}\cos\left(2\alpha x_2\right)+\frac{\alpha\sqrt{3}}{16}\cos\left(2\alpha x_3\right)+\nonumber\\
&-&\frac{3\alpha}{32}\sin\left(2\alpha x_2\right)+\frac{\alpha\sqrt{3}}{32}\cos\left(2\alpha x_2\right)-\nonumber\\
&-&\frac{3\alpha}{32}\sin\left(2\alpha x_3\right)-\frac{\alpha\sqrt{3}}{32}\cos\left(2\alpha x_3\right)\nonumber\\
&=&0
\end{eqnarray}
%
\begin{eqnarray}
\label{4.20}
\frac{\mathcal{V}_2^0\left( {x_3,x_1;-\frac{\pi}{3}, \frac{\pi}{3}, \frac{\pi}{2}} \right)}{v_r^2}=0
\end{eqnarray}
%
\begin{eqnarray}
\label{4.21}
\frac{\mathcal{V}_3^0\left( {x_1,x_2;-\frac{\pi}{3}, \frac{\pi}{3}, \frac{\pi}{2}} \right)}{v_r^2}=0
\end{eqnarray}
In light of the fact that $g_{i}^0$ equals $\mathcal{W}_i^0$, the following expressions are derived for $g_i^0$ and $\mathcal{G}^0\left(x_1, x_2,x_3\right)=\sum_{k = 1}^3 {\frac{{\partial \mathcal{W}_k^0 }}{{\partial x_k }}}$:
%
\begin{eqnarray}
\label{4.22}
\frac{g_{1}^0\left( {x_1,x_2,x_3} \right)}{{v_r^2}}&\!=\!&-\frac{3\alpha }{8}\sin\left(2\alpha x_1 \right)\!-\!\frac{3\alpha }{32}\sin \left(2\alpha x_1 \right)\cos \left(2\alpha x_2 \right)\!-\!\frac{3\alpha }{32}\sin \left(2\alpha x_1 \right)\cos \left(2\alpha x_3 \right)-\nonumber\\
&-&\frac{3\sqrt{3}\alpha }{32}\cos \left(2\alpha x_1 \right)\cos \left(2\alpha x_3 \right)+\frac{3\sqrt{3} \alpha }{32}\cos \left(2\alpha x_1 \right)\cos \left(2\alpha x_2 \right)+\nonumber\\
&+&\frac{3\sqrt{3} \alpha }{32}\sin \left(2\alpha x_1 \right)\sin \left(2\alpha x_2 \right)-\frac{3\sqrt{3} \alpha}{32}\sin \left(2\alpha x_1 \right)\sin \left(2\alpha x_3 \right)+\nonumber\\
&+&\frac{3\alpha }{32}\cos \left(2\alpha x_1 \right)\sin \left(2\alpha x_3 \right)+\frac{3\alpha }{32}\cos \left(2\alpha x_1 \right)\sin \left(2\alpha x_2 \right)
\end{eqnarray}
%
\begin{eqnarray}
\label{4.23}
g_{2}^0\left( {x_2,x_3,x_1} \right)&=&v_r^2g_{1}^0\left( {x_2,x_3,x_1} \right)
\end{eqnarray}
%
\begin{eqnarray}
\label{4.24}
g_{3}^0\left( {x_3,x_1,x_2} \right)&=&v_r^2g_{2}^0\left( {x_3,x_1,x_2} \right)
\end{eqnarray}
%
\begin{eqnarray}
\label{4.25}
\frac{\mathcal{G}^0\left(x_1, x_2,x_3\right)}{{v_r^2}}&=&-\frac{3\alpha^2 }{4}\cos\left(2\alpha x_1 \right)-\frac{3\alpha^2 }{4}\cos\left(2\alpha x_2 \right)-\frac{3\alpha^2 }{4}\cos\left(2\alpha x_3 \right)-\nonumber\\
&-&\frac{3\alpha^2 }{8}\cos \left(2\alpha x_1 \right)\cos \left(2\alpha x_2 \right)-\frac{3\alpha^2 }{8}\cos \left(2\alpha x_1 \right)\cos \left(2\alpha x_3 \right)+\nonumber\\
&+&\frac{3\sqrt{3}\alpha^2 }{8}\sin \left(2\alpha x_1 \right)\cos \left(2\alpha x_3 \right)-\frac{3\sqrt{3} \alpha^2 }{8}\sin \left(2\alpha x_1 \right)\cos \left(2\alpha x_2 \right)+\nonumber\\
&+&\frac{3\sqrt{3} \alpha^2 }{8}\cos \left(2\alpha x_1 \right)\sin \left(2\alpha x_2 \right)-\frac{3\sqrt{3} \alpha^2}{8}\cos \left(2\alpha x_1 \right)\sin \left(2\alpha x_3 \right)-\nonumber\\
&-&\frac{3\alpha^2 }{8}\sin \left(2\alpha x_1 \right)\sin \left(2\alpha x_3 \right)-\frac{3\alpha^2 }{8}\sin \left(2\alpha x_1 \right)\sin \left(2\alpha x_2 \right)-\nonumber\\
&-&\frac{3\alpha^2 }{8}\cos \left(2\alpha x_2 \right)\cos \left(2\alpha x_3 \right)-\frac{3\sqrt{3} \alpha^2 }{8}\sin \left(2\alpha x_2 \right)\cos \left(2\alpha x_3 \right)+\nonumber\\
&+&\frac{3\sqrt{3} \alpha^2 }{8}\cos \left(2\alpha x_2 \right)\sin \left(2\alpha x_3 \right)-\frac{3\alpha^2 }{8}\sin \left(2\alpha x_2 \right)\sin \left(2\alpha x_3 \right)
\end{eqnarray}
By substituting $\mathcal{G}^0\left(x_1, x_2,x_3\right)$ from $\left(\ref{4.25}\right)$ into  $\left(\ref{4.5}\right)$ and evaluating the definite integrals term by term, we can see that the results equal $g_i^0$ derived from $\left(\ref{2.6}\right)$; that is,  $\mathcal{U}_i^0 \equiv 0$, $i = 1,2,3$. We have used the integral identities $\left(\ref{B4}\right)$-$\left(\ref{B8}\right)$ given in Appendix B to perform the integrations.\\\\ 
Substituting for the initial condition $v_i^{0}$, $i = 1,2,3$ from $\left(\ref{4.1}\right)$, $\left(\ref{4.2}\right)$ and $\left(\ref{4.3}\right)$ into $\left(\ref{3.12}\right)$ with phase angles $\xi_1$, $\xi_2$ and $\xi_3 $ set, respectively, to $-\frac{\pi}{3}$, $\frac{\pi}{3}$, and $\frac{\pi}{2}$ and performing the integrations, we arrive at the solution of the Navier-Stokes equation $\left(\ref{2.1}\right)$ in three dimensions:
%
\begin{eqnarray}
\label{4.26}
v_1& = &{v_r}\left[\sin { \left(\alpha x_1-\frac{\pi}{3} \right) } \cos {\left(\alpha x_2+\frac{\pi}{3}\right) }\sin { \left(\alpha x_3+\frac{\pi}{2}\right) }\right. -\nonumber\\
&-&\left.\sin  {\left( \alpha x_1+\frac{\pi}{3}\right) } \cos  { \left(\alpha x_3-\frac{\pi}{3} \right) } \sin  {\left( \alpha x_2+\frac{\pi}{2}\right) }\right]e^{ - 3\alpha ^2 \kappa t}
\end{eqnarray}
%
\begin{eqnarray}
\label{4.27}
v_2 & = &{v_r}\left[\sin { \left(\alpha x_2-\frac{\pi}{3} \right) } \cos {\left(\alpha x_3+\frac{\pi}{3}\right) }\sin { \left(\alpha x_1+\frac{\pi}{2}\right) }\right.-\nonumber\\
&-&\left.\sin  {\left( \alpha x_2+\frac{\pi}{3}\right) } \cos  { \left(\alpha x_1-\frac{\pi}{3} \right) } \sin  {\left( \alpha x_3+\frac{\pi}{2}\right) }\right]e^{ - 3\alpha ^2 \kappa t}
\end{eqnarray}
%
\begin{eqnarray}
\label{4.28}
v_3 &=&{v_r}\left[ \sin { \left(\alpha x_3-\frac{\pi}{3} \right) } \cos {\left(\alpha x_1+\frac{\pi}{3}\right) }\sin { \left(\alpha x_2+\frac{\pi}{2}\right) }\right. -\nonumber\\
&-&\left.\sin  {\left( \alpha x_3+\frac{\pi}{3}\right) } \cos  { \left(\alpha x_2-\frac{\pi}{3} \right) } \sin  {\left( \alpha x_1+\frac{\pi}{2}\right) }\right]e^{ - 3\alpha ^2 \kappa t} 
\end{eqnarray}
We have used the integral identities $\left(\ref{B1}\right)$-$\left(\ref{B3}\right)$ given in Appendix B to perform the integrations.\\\\ 
Pressure is obtained from the solution of the Poisson equation $\left(\ref{3.2}\right)$. Because $\left(\ref{3.10}\right)$ is satisfied, we can express $\left(\ref{3.5}\right)$ as follows:
%
\begin{eqnarray}
\label{4.29}
\frac{{\partial p }}{{\partial x_i }} = -\rho g_i=-\rho g_i^0e^{ - 6\alpha ^2 \kappa t} \qquad i = 1,2,3
\end{eqnarray}
This leads to straightforward integration for pressure:
%
\begin{eqnarray}
\label{4.30}
p&=&-\frac{3\rho {v_r^2} }{16}\left[\cos\left(2\alpha x_1 \right)+\cos\left(2\alpha x_2 \right)+\cos\left(2\alpha x_3 \right)\right]e^{ - 6\alpha ^2 \kappa t}-\nonumber\\
&-&\frac{3\rho{v_r^2}}{64}\left[\cos \left(2\alpha x_1-2\alpha x_2 \right)+\cos \left(2\alpha x_3-2\alpha x_1 \right)+\cos \left(2\alpha x_3-2\alpha x_2 \right)\right]e^{ - 6\alpha ^2 \kappa t}-\nonumber\\
&-&\frac{3\sqrt{3}\rho{v_r^2} }{64}\left[\sin \left(2\alpha x_1-2\alpha x_2 \right)+\sin \left(2\alpha x_3-2\alpha x_1 \right)+\sin \left(2\alpha x_2-2\alpha x_3\right)\right]e^{ - 6\alpha ^2 \kappa t}\qquad\quad
\end{eqnarray}
\section{A solution in  $\mathbb{R}^4=\left\{ { - \infty  < x_i  < \infty ; \,\, i = 1,2,3,4} \right\}$}
At time zero, we consider the solenoidal velocity vector field $v_i^{0}\left( x \right)$ of a quadruple-periodic form 
\begin{eqnarray}
\label{5.1}
\frac{v_1^0}{v_r} &=& \sin { \left(\alpha x_1+\xi_1 \right) } \cos {\left(\alpha x_2+\xi_2\right) }\sin { \left(\alpha x_3+\xi_3\right) } \cos { \left(\alpha x_4+\xi_4 \right) } -\nonumber\\
&-&\sin  {\left( \alpha x_1+\xi_2\right) }\sin  {\left( \alpha x_2+\xi_3\right) } \cos  { \left(\alpha x_3+\xi_4 \right) } \cos  {\left( \alpha x_4+\xi_1\right) } 
\end{eqnarray}
%
\begin{eqnarray}
\label{5.2}
\frac{v_2^0}{v_r} &=& \sin { \left(\alpha x_2+\xi_1 \right) } \cos {\left(\alpha x_3+\xi_2\right) }\sin { \left(\alpha x_4+\xi_3\right) } \cos { \left(\alpha x_1+\xi_4 \right) } -\nonumber\\
&-&\sin  {\left( \alpha x_2+\xi_2\right) }\sin  {\left( \alpha x_3+\xi_3\right) } \cos  { \left(\alpha x_4+\xi_4 \right) } \cos  {\left( \alpha x_1+\xi_1\right) } 
\end{eqnarray}
%
\begin{eqnarray}
\label{5.3}
\frac{v_3^0}{v_r}&=& \sin { \left(\alpha x_3+\xi_1 \right) } \cos {\left(\alpha x_4+\xi_2\right) }\sin { \left(\alpha x_1+\xi_3\right) } \cos { \left(\alpha x_2+\xi_4 \right) }  -\nonumber\\
&-&\sin  {\left( \alpha x_3+\xi_2\right) } \sin  {\left( \alpha x_4+\xi_3\right) } \cos  { \left(\alpha x_1+\xi_4 \right) }\cos  {\left( \alpha x_2+\xi_1\right) }  
\end{eqnarray}
%
\begin{eqnarray}
\label{5.4}
\frac{v_4^0}{v_r} &=& \sin { \left(\alpha x_4+\xi_1 \right) } \cos {\left(\alpha x_1+\xi_2\right) }\sin { \left(\alpha x_2+\xi_3\right) } \cos { \left(\alpha x_3+\xi_4 \right) } -\nonumber\\
&-&\sin  {\left( \alpha x_4+\xi_2\right) }\sin  {\left( \alpha x_1+\xi_3\right) } \cos  {\left( \alpha x_3+\xi_1\right) } \cos  { \left(\alpha x_2+\xi_4 \right) } 
\end{eqnarray}
where $\xi_1$, $\xi_2$, $\xi_3$ and $\xi_4$ are phase angles. The the externally applied force on the fluid $f_i$ is set to zero. Equation $\left(\ref{3.10}\right)$ may be written as follows at time zero in $\mathbb{R}^4$:
\begin{eqnarray}
\label{5.5}
g_i^{0}= \frac{1}{{2\pi^2 }}\int\limits_{\mathbb{R}^4 } {\frac{{\left( {x_i  - y_i } \right)\mathcal{G}^0\left(y_1, y_2,y_3,y_4\right) }}
{{\left\{ {\mathcal{P}_4\left( {x,y} \right)} \right\}^2 }}} \prod\limits_{j = 1}^4 {dy_j }\qquad i = 1,2,3,4
\end{eqnarray}
where $\mathcal{G}^0\left(y_1, y_2,y_3,y_4\right)=\sum_{k = 1}^4 \frac{\partial g_k^0 \left( {y,t} \right)}{\partial y_k }$. Changing the variable of integration, $u_i=x_i  - y_i$, gives
%
\begin{eqnarray}
\label{5.6}
g_i^{0}= \frac{1}{{2\pi^2 }}\int\limits_{\mathbb{R}^4 } {\frac{{u_i\,\mathcal{G}^0\left(x_1-u_1, x_2-u_2,x_3-u_3,x_4-u_4\right) }}
{{\left[u_1^2+u_2^2+u_3^2+u_4^2 \right]^2 }}} \prod\limits_{j = 1}^4 {du_j }\nonumber\\i = 1,2,3,4
\end{eqnarray}
Following the method of the previous case, we divide $g_i^0$ into two parts:
%
\begin{eqnarray}
\label{5.7}
g_i^0 = \mathcal{V}_i^0 + \mathcal{W}_i^0\qquad i = 1,2,3,4
\end{eqnarray}
By substituting $v_1^0$, $v_2^0$, $v_3^0$ and $v_4^0$ into $\left(\ref{2.6}\right)$ and separating the terms into  $\mathcal{V}_i^0$ and $\mathcal{W}_i^0$ we get the following expressions:
%
\begin{eqnarray}
\label{5.8}
g_{1}^0\left( {x_1,x_2,x_3,x_4; \xi_1, \xi_2, \xi_3, \xi_4} \right)&=&\mathcal{V}_1^0\left( {x_2,x_3,x_4; \xi_1, \xi_2, \xi_3, \xi_4} \right) +\nonumber\\
&+& \mathcal{W}_1^0\left( {x_1,x_2,x_3,x_4; \xi_1, \xi_2, \xi_3, \xi_4} \right)
\end{eqnarray}
%
\begin{eqnarray}
\label{5.9}
g_{2}^0\left( {x_2,x_3,x_4,x_1; \xi_1, \xi_2, \xi_3, \xi_4} \right)&=&\mathcal{V}_2^0\left( {x_3,x_4,x_1; \xi_1, \xi_2, \xi_3, \xi_4} \right) +\nonumber\\
&+& \mathcal{W}_2^0\left( {x_2,x_3,x_4,x_1; \xi_1, \xi_2, \xi_3, \xi_4} \right)
\end{eqnarray}
%
\begin{eqnarray}
\label{5.10}
g_{3}^0\left( {x_3,x_4,x_1,x_2; \xi_1, \xi_2, \xi_3, \xi_4} \right)&=&\mathcal{V}_3^0\left( {x_4,x_1,x_2; \xi_1, \xi_2, \xi_3, \xi_4} \right) + \nonumber\\
&+&\mathcal{W}_3^0\left( {x_3,x_4,x_1,x_2; \xi_1, \xi_2, \xi_3, \xi_4} \right)
\end{eqnarray}
%
\begin{eqnarray}
\label{5.11}
g_{4}^0\left( {x_4,x_1,x_2,x_3; \xi_1, \xi_2, \xi_3, \xi_4} \right)&=&\mathcal{V}_4^0\left( {x_1,x_2,x_3; \xi_1, \xi_2, \xi_3, \xi_4} \right) +\nonumber\\
&+& \mathcal{W}_4^0\left( {x_4,x_1,x_2,x_3; \xi_1, \xi_2, \xi_3, \xi_4} \right)
\end{eqnarray}
The expressions for $\frac{\mathcal{V}_1^0\left( {x_2,x_3,x_4; \xi_1, \xi_2, \xi_3, \xi_4} \right)}{v_r^2}$ and $\frac{\mathcal{W}_1^0\left( {x_1,x_2,x_3,x_4; \xi_1, \xi_2, \xi_3, \xi_4} \right)}{v_r^2}$ are given by $\left(\ref{C1}\right)$ and $\left(\ref{C2}\right)$ respectively in Appendix C.
%
\begin{eqnarray}
\label{5.12}
\mathcal{V}_2^0\left( {x_3,x_4,x_1; \xi_1, \xi_2, \xi_3, \xi_4} \right)&=&v_r^2\mathcal{V}_1^0\left( {x_3,x_4,x_1; \xi_1, \xi_2, \xi_3, \xi_4} \right)
\end{eqnarray}
%
\begin{eqnarray}
\label{5.13}
\mathcal{V}_3^0\left( {x_4,x_1,x_2; \xi_1, \xi_2, \xi_3, \xi_4} \right)&=&v_r^2\mathcal{V}_2^0\left( {x_4,x_1,x_2; \xi_1, \xi_2, \xi_3, \xi_4} \right)
\end{eqnarray}
%
\begin{eqnarray}
\label{5.14}
\mathcal{V}_4^0\left( {x_1,x_2,x_3; \xi_1, \xi_2, \xi_3, \xi_4} \right)&=&v_r^2\mathcal{V}_3^0\left( {x_1,x_2,x_3; \xi_1, \xi_2, \xi_3, \xi_4} \right) 
\end{eqnarray}
%
\begin{eqnarray}
\label{5.15}
\mathcal{W}_2^0\left( {x_2,x_3,x_4,x_1; \xi_1, \xi_2, \xi_3, \xi_4} \right)&=&v_r^2\mathcal{W}_1^0\left( {x_2,x_3,x_4,x_1; \xi_1, \xi_2, \xi_3, \xi_4} \right)
\end{eqnarray}
%
\begin{eqnarray}
\label{5.16}
\mathcal{W}_3^0\left( {x_3,x_4,x_1,x_2; \xi_1, \xi_2, \xi_3, \xi_4} \right)&=&v_r^2\mathcal{W}_2^0\left( {x_3,x_4,x_1,x_2; \xi_1, \xi_2, \xi_3, \xi_4} \right)
\end{eqnarray}
%
\begin{eqnarray}
\label{5.17}
\mathcal{W}_4^0\left( {x_4,x_1,x_2,x_3; \xi_1, \xi_2, \xi_3, \xi_4} \right)&=&v_r^2\mathcal{W}_3^0\left( {x_4,x_1,x_2,x_3; \xi_1, \xi_2, \xi_3, \xi_4} \right)
\end{eqnarray}
If $\left(\ref{5.6}\right)$ is to hold the following prerequisites must be met:
\begin{eqnarray}
\label{5.18}
\mathcal{V}_1^0\left( {x_2,x_3,x_4; \xi_1, \xi_2, \xi_3, \xi_4} \right) =0,\,\, \forall -\infty<x_2<\infty,\,\,-\infty<x_3<\infty\,\,\textnormal{and}\,\,-\infty<x_4<\infty\nonumber\\
\end{eqnarray}
%
\begin{eqnarray}
\label{5.19}
\mathcal{V}_2^0\left( {x_3,x_4,x_1; \xi_1, \xi_2, \xi_3, \xi_4} \right)=0,\,\, \forall -\infty<x_3<\infty,\,\,-\infty<x_4<\infty\,\,\textnormal{and}\,\,-\infty<x_1<\infty\nonumber\\
\end{eqnarray}
%
\begin{eqnarray}
\label{5.20}
\mathcal{V}_3^0\left( {x_4,x_1,x_2; \xi_1, \xi_2, \xi_3, \xi_4} \right)=0,\,\, \forall -\infty<x_4<\infty,\,\,-\infty<x_1<\infty\,\,\textnormal{and}\,\,-\infty<x_2<\infty\nonumber\\
\end{eqnarray}
%
\begin{eqnarray}
\label{5.21}
\mathcal{V}_4^0\left( {x_1,x_2,x_3; \xi_1, \xi_2, \xi_3, \xi_4} \right)=0,\,\, \forall -\infty<x_1<\infty,\,\,-\infty<x_2<\infty\,\,\textnormal{and}\,\,-\infty<x_3<\infty\nonumber\\
\end{eqnarray}
The phase angles $\xi_1$, $\xi_2$, $\xi_3$ and $\xi_4$ are chosen such that they satisfy the prerequisites $\left(\ref{5.18}\right)$, $\left(\ref{5.19}\right)$, $\left(\ref{5.20}\right)$ and $\left(\ref{5.21}\right)$. In $\left(\ref{C1}\right)$, $\left(\ref{5.12}\right)$, $\left(\ref{5.13}\right)$ and $\left(\ref{5.14}\right)$, the phase angles $\xi_1=-\frac{\pi}{4}$, $\xi_2=\frac{\pi}{4}$, $\xi_3=\frac{\pi}{4}$ and $\xi_4=\frac{\pi}{4}$ satisfy  $\left(\ref{5.18}\right)$, $\left(\ref{5.19}\right)$, $\left(\ref{5.20}\right)$ and $\left(\ref{5.21}\right)$:
%
\begin{eqnarray}
\label{5.22}
\frac{\mathcal{V}_1^0\left( {x_2,x_3,x_4; -\frac{\pi}{4}, \frac{\pi}{4}, \frac{\pi}{4}, \frac{\pi}{4}} \right)}{v_r^2}=0
\end{eqnarray}
%
\begin{eqnarray}
\label{5.23}
\frac{\mathcal{V}_2^0\left( {x_3,x_4,x_1;  -\frac{\pi}{4}, \frac{\pi}{4}, \frac{\pi}{4}, \frac{\pi}{4}} \right)}{v_r^2}=0
\end{eqnarray}
%
\begin{eqnarray}
\label{5.24}
\frac{\mathcal{V}_3^0\left( {x_4,x_1,x_2;  -\frac{\pi}{4}, \frac{\pi}{4}, \frac{\pi}{4}, \frac{\pi}{4}} \right)}{v_r^2}=0
\end{eqnarray}
%
\begin{eqnarray}
\label{5.25}
\frac{\mathcal{V}_4^0\left( {x_1,x_2,x_3;   -\frac{\pi}{4}, \frac{\pi}{4}, \frac{\pi}{4}, \frac{\pi}{4}}  \right)}{v_r^2}=0
\end{eqnarray}
Taking note that $g_{i}^0$ equals $\mathcal{W}_i^0$, the following expressions are derived for $g_i^0$ and\\ $\mathcal{G}^0\left(x_1, x_2,x_3,x_3\right)=\sum_{k = 1}^4 {\frac{{\partial \mathcal{W}_k^0 }}{{\partial x_k }}}$:
%
\begin{eqnarray}
\label{5.26}
g_1^0  &=& -\frac{v_r^2\alpha}{2}\cos  { \left(2\alpha x_1 \right) } \sin  {\left( 2\alpha x_3\right) }
\end{eqnarray}
%
\begin{eqnarray}
\label{5.27}
g_2^0  &=& -\frac{v_r^2\alpha}{2}\cos  { \left(2\alpha x_2 \right) } \sin  {\left( 2\alpha x_4\right) }
\end{eqnarray}
%
\begin{eqnarray}
\label{28}
g_3^0  &=& -\frac{v_r^2\alpha}{2}\cos  { \left(2\alpha x_3 \right) } \sin  {\left( 2\alpha x_1\right) }
\end{eqnarray}
%
\begin{eqnarray}
\label{22}
g_4^0  &=& -\frac{v_r^2\alpha}{2}\cos  { \left(2\alpha x_4 \right) } \sin  {\left( 2\alpha x_2\right) }
\end{eqnarray}
%
\begin{eqnarray}
\label{5.30}
\frac{\mathcal{G}^0\left(x_1, x_2,x_3,x_4\right)}{{v_r^2}}&=&\alpha^2\sin  { \left(2\alpha x_1 \right) } \sin  {\left( 2\alpha x_3\right) }+\alpha^2\sin  { \left(2\alpha x_2 \right) } \sin  {\left( 2\alpha x_4\right) }+ \nonumber\\
&+&\alpha^2\sin  { \left(2\alpha x_3 \right) } \sin  {\left( 2\alpha x_1\right) }+\alpha^2\sin  { \left(2\alpha x_4 \right) } \sin  {\left( 2\alpha x_2\right) }
\end{eqnarray}
By substituting $\mathcal{G}^0\left(x_1, x_2,x_3,x_4\right)$ from $\left(\ref{5.30}\right)$ into  $\left(\ref{5.6}\right)$ and evaluating the definite integrals term by term, we can see that the results equal $g_i^0$ derived from $\left(\ref{2.6}\right)$; that is,  $\mathcal{U}_i^0 \equiv 0$, $i = 1,2,3,4$.. We have used the integral identity $\left(\ref{B9}\right)$ in Appendix B to perform the integrations.\\\\ 
Substituting for the initial condition $v_i^{0}$, $i = 1,2,3,4$ from $\left(\ref{5.1}\right)$, $\left(\ref{5.2}\right)$, $\left(\ref{5.3}\right)$ and $\left(\ref{5.4}\right)$ into $\left(\ref{3.12}\right)$ with phase angles $\xi_1$, $\xi_2$, $\xi_3 $ and $\xi_4$ set, respectively, to $\xi_1=-\frac{\pi}{4}$, $\xi_2=\frac{\pi}{4}$, $\xi_3=\frac{\pi}{4}$ and $\xi_4=\frac{\pi}{4}$ and performing the integrations, we arrive at the solution of the Navier-Stokes equation $\left(\ref{2.1}\right)$ in four dimensions:
%
\begin{eqnarray}
\label{5.31}
v_1 &=& v_r\left[\sin { \left(\alpha x_1-\frac{\pi}{4} \right) } \cos {\left(\alpha x_2+\frac{\pi}{4}\right) }\sin { \left(\alpha x_3+\frac{\pi}{4}\right) } \cos { \left(\alpha x_4+\frac{\pi}{4} \right) }-\right.\nonumber\\
&-&\left.\sin  {\left( \alpha x_1+\frac{\pi}{4}\right) }\sin  {\left( \alpha x_2+\frac{\pi}{4}\right) } \cos  { \left(\alpha x_3+\frac{\pi}{4} \right) } \cos  {\left( \alpha x_4-\frac{\pi}{4}\right) } \right]e^{ - 4\alpha ^2 \kappa t}\qquad
\end{eqnarray}
%
\begin{eqnarray}
\label{5.32}
v_2 &=&v_r \left[\sin { \left(\alpha x_2-\frac{\pi}{4} \right) } \cos {\left(\alpha x_3+\frac{\pi}{4}\right) }\sin { \left(\alpha x_4+\frac{\pi}{4}\right) } \cos { \left(\alpha x_1+\frac{\pi}{4} \right) }-\right.\nonumber\\
&-&\left.\sin  {\left( \alpha x_2+\frac{\pi}{4}\right) }\sin  {\left( \alpha x_3+\frac{\pi}{4}\right) } \cos  { \left(\alpha x_4+\frac{\pi}{4} \right) } \cos  {\left( \alpha x_1-\frac{\pi}{4}\right) } \right]e^{ - 4\alpha ^2 \kappa t}\qquad
\end{eqnarray}
%
\begin{eqnarray}
\label{5.33}
v_3&=&v_r\left[ \sin { \left(\alpha x_3-\frac{\pi}{4} \right) } \cos {\left(\alpha x_4+\frac{\pi}{4}\right) }\sin { \left(\alpha x_1+\frac{\pi}{4}\right) } \cos { \left(\alpha x_2+\frac{\pi}{4} \right) } -\right.\nonumber\\
&-&\left.\sin  {\left( \alpha x_3+\frac{\pi}{4}\right) } \sin  {\left( \alpha x_4+\frac{\pi}{4}\right) } \cos  { \left(\alpha x_1+\frac{\pi}{4} \right) }\cos  {\left( \alpha x_2-\frac{\pi}{4}\right) }\right]e^{ - 4\alpha ^2 \kappa t}\qquad 
\end{eqnarray}
%
\begin{eqnarray}
\label{5.34}
v_4&=&v_r\left[ \sin { \left(\alpha x_4-\frac{\pi}{4} \right) } \cos {\left(\alpha x_1+\frac{\pi}{4}\right) }\sin { \left(\alpha x_2+\frac{\pi}{4}\right) } \cos { \left(\alpha x_3+\frac{\pi}{4} \right) }-\right.\nonumber\\
&-&\left.\sin  {\left( \alpha x_4+\frac{\pi}{4}\right) }\sin  {\left( \alpha x_1+\frac{\pi}{4}\right) }  \cos  { \left(\alpha x_2+\frac{\pi}{4} \right) }\cos  {\left( \alpha x_3-\frac{\pi}{4}\right) } \right]e^{ - 4\alpha ^2 \kappa t}\qquad
\end{eqnarray}
We have used the integral identities $\left(\ref{B1}\right)$-$\left(\ref{B3}\right)$ given in Appendix B to perform the integrations.\\\\ 
Pressure is obtained from the solution of the Poisson equation $\left(\ref{3.2}\right)$. Because $\left(\ref{3.10}\right)$ is satisfied, we can express $\left(\ref{3.5}\right)$ as follows:
%
\begin{eqnarray}
\label{5.35}
\frac{{\partial p }}{{\partial x_i }} =- \rho g_i=-\rho g_i^0e^{ - 8\alpha ^2 \kappa t} \qquad i = 1,2,3,4
\end{eqnarray}
This leads to straightforward integration for pressure:
%
\begin{eqnarray}
\label{5.36}
p&=&\frac{\rho {v_r^2} }{4}\left[\sin  { \left(2\alpha x_1 \right) } \sin  {\left( 2\alpha x_3\right) }+\sin  { \left(2\alpha x_2 \right) } \sin  {\left( 2\alpha x_4\right) }\right]e^{ - 8\pi ^2 \kappa t}
\end{eqnarray}
%
\vspace*{-0.50 in}
\section{Concluding Remarks} 
The statement of the problem and the edifice of the proposed solution method in a space having an arbitrary number of dimensions,  $\mathbb{R}^n$, are presented in Sections 2 and 3, respectively. The pressure field is formally obtained by taking the divergence of the Navier-Stokes equation $\left(\ref{2.1}\right)$ as a solution to the Poisson equation, yielding an expression $\left(\ref{3.5}\right)$ for pressure gradient. By inserting the pressure gradient back into $\left(\ref{2.1}\right)$ and equating the unsteady terms, $\frac{{\partial v_i }}{{\partial t}}$, to the sum of  three terms that are associated, respectively, with the linear viscous force,  the nonlinear inertial force and the externally applied force acting on the fluid we obtain the rephrased Navier-Stokes equation $\left(\ref{3.7}\right)$.  Subject to satisfying an integral equation $\left(\ref{3.10}\right)$, $\left(\ref{3.7}\right)$ is then reduced to a non-homogeneous diffusion equation in velocity. When the externally applied force is set to zero, it is further reduced to the Cauchy diffusion equation. The nonlinearity is entirely isolated to the integral equation $\left(\ref{3.10}\right)$.\\\\ 
Succinctly put, subject to satisfying a certain condition $\left(\ref{3.10}\right)$, the Navier-Stokes system of equations is decomposed into two linear equations, the non-homogeneous diffusion equation and the Poisson equation, for the velocity and pressure fields respectively.  The velocity and pressure fields are given by the solutions of Cauchy diffusion equation and Poisson equation respectively. The proposed solution methods belong to the class of Beltrami flows and are valid under proper regularity conditions at infinity. \\\\ 
The technique presented in this paper is a straightforward and direct method to derive solutions of the velocity and pressure fields for the Navier-Stokes flow problem. We have demonstrated its efficacy by applying it to reproduce the double-periodic solution of \cite{tay}, the unit-periodic solution of \cite{arn}, and the triple-periodic solution of \cite{ant}. In addition, we have also derived a new quadruple-periodic form of the solution in four dimensions. The method is easy to master and can almost be reduced to a \enquote*{drill}.
\appendix
%
%
\section{Derivations of two known solutions \cite{tay, arn} satisfying equation $\left(\ref{3.10}\right)$ presented from \cite{tha1}} 
\subsection{The Taylor solution in $\mathbb{R}^2=\left\{ { - \infty  < x_i  < \infty ;  i = 1,2} \right\}$}
A two-dimensional solenoidal velocity vector field of a double-periodic form $v_i^{0}\left( x \right)$ is considered at time zero: 
%
\begin{eqnarray}
\label{A1}
v_1^{0}  =v_r \sin \left( {\pi x_1 } \right)\cos \left( {\pi x_2 } \right)
\end{eqnarray}
%
\begin{eqnarray}
\label{A2}
v_2^{0}  =  -v_r \cos \left( {\pi x_1 } \right)\sin \left( {\pi x_2 } \right)
\end{eqnarray}
The externally applied force on the fluid $f_i$ is set to zero. At time zero in $\mathbb{R}^2$,  $\left(\ref{3.10}\right)$ may be written as follows:
%
\begin{eqnarray}
\label{A3}
g_i^{0}=\frac{1}
{{2\pi }}\int\limits_{\mathbb{R}^2 } {\frac{{\left( {x_i  - y_i } \right)\sum\limits_{k = 1}^2 {\mathcal{G}^0\left(y_1, y_2\right)} }}
{{\mathcal{P}_2 \left( {x,y} \right)}}} \prod\limits_{j = 1}^2 {dy_j } \qquad i = 1,2
\end{eqnarray}
Where $g_i^0  = \sum_{j = 1}^2 {v_j^0 \frac{{\partial v_i^0}}{{\partial x_j }}}$ and $\mathcal{G}^0\left(y_1, y_2\right)=\sum_{k = 1}^2 \frac{\partial g_k^0 \left( {y,t} \right)}{\partial y_k }$. Changing the variable of integration, $u_i=x_i  - y_i$ and performing the integrations we obtain 
%
\begin{eqnarray}
\label{A4}
 g_i^{0}= \frac{1}
{{2\pi }}\int\limits_{\mathbb{R}^2 }{\frac{{u_i\,\mathcal{G}^0\left(x_1-u_1, x_2-u_2\right) }}
{{\left[u_1^2+u_2^2\right] }}} \prod\limits_{j = 1}^2 {du_j }= \frac{\pi }
{2}\sin \left( {2\pi x_i } \right) \quad i = 1,2
\end{eqnarray}
resulting in $\mathcal{U}_i\left( x,t \right) \equiv 0$.  $v\left( {x,t} \right)$ is obtained from $\left(\ref{3.12}\right)$:
%
\begin{eqnarray}
\label{A5}
\frac{v_1}{{v_r}}  = \sin \left( {\pi x_1 } \right)\cos \left( {\pi x_2 } \right)e^{ - 2\pi ^2 \kappa t}
\end{eqnarray}
%
\begin{eqnarray}
\label{A6}
\frac{v_2}{{v_r}}  =  - \cos \left( {\pi x_1 } \right)\sin \left( {\pi x_2 } \right)e^{ - 2\pi ^2 \kappa t}
\end{eqnarray}
We have used the integral identities $\left(\ref{B1}\right)$-$\left(\ref{B3}\right)$ given in Appendix B to perform the integrations. The pressure $p\left( {x,t} \right)$ is obtained by integrating (3.13):
\begin{eqnarray}
\label{A7}
p =  \frac{{\rho {v_r^2}e^{ - 4\pi ^2 \kappa t} }}{4}\left[ {\cos \left( {2\pi x_1 } \right) + \cos \left( {2\pi x_2 } \right)} \right]
\end{eqnarray}
\subsection{Arnold-Beltrami-Childress (ABC) flows in $\mathbb{R}^3=\left\{ { - \infty  < x_i  < \infty ; \,\, i = 1,2,3} \right\}$}
A three-dimensional solenoidal velocity vector field of a unit-periodic form $v_i^{0}\left( x \right)$ is considered at time zero: 
%
\begin{eqnarray}
\label{A8}
\frac{v_1^{(0)}}{v_r}  = a\sin \pi x_3  - c\cos \pi x_2
\end{eqnarray}
%
\begin{eqnarray}
\label{A9}
\frac{v_2^{(0)}}{v_r}   = b\sin \pi x_1 - a\cos \pi x_3
\end{eqnarray}
%
\begin{eqnarray}
\label{A10}
\frac{v_3^{(0)}}{v_r}  = c\sin \pi x_2  - b\cos \pi x_1
\end{eqnarray}
Where $a$, $b$ and $c$ are real constants. The externally applied force on the fluid $f_i$ is set to zero. At time zero in $\mathbb{R}^3$,  $\left(\ref{3.10}\right)$ may be written as follows
%
\begin{eqnarray}
\label{A11}
g_i^{0}=  \frac{1}
{{4\pi }}\int\limits_{\mathbb{R}^3 } {\frac{{\left( {x_i  - y_i } \right)\mathcal{G}^0\left(y_1, y_2,y_3\right) }}
{{\left\{ {\mathcal{P}_3\left( {x,y} \right)} \right\}^{\frac{3}
{2}} }}} \prod\limits_{j = 1}^3 {dy_j }\qquad i = 1,2,3
\end{eqnarray}
Where $g_i^0  = \sum_{j = 1}^3 {v_j^0 \frac{{\partial v_i^0}}{{\partial x_j }}}$ and $\mathcal{G}^0\left(y_1, y_2\right)=\sum_{k = 1}^3 \frac{\partial g_k^0 \left( {y,t} \right)}{\partial y_k }$. Changing the variable of integration, $u_i=x_i  - y_i$ and performing the integrations we obtain 
%
\begin{eqnarray}
\label{A12}
g_1^{(0)}  &=&  \frac{1}{{4\pi }}\int\limits_{\mathbb{R}^3 }{\frac{{u_1\,\mathcal{G}^0\left(x_1-u_1, x_2-u_2,x_3-u_3\right) }}
{{\left[u_1^2+u_2^2+u_3^2 \right]^{\frac{3}{2}} }}}  \prod\limits_{j = 1}^3 {du_j }
\nonumber\\
&  = & 
\pi \left\{ {bc\sin \left(\pi x_1\right) \sin \left(\pi x_2\right)- ab\cos \left(\pi x_1\right) \cos\left( \pi x_3\right) } \right\}\qquad
\end{eqnarray}
%
\begin{eqnarray}
\label{A13}
g_2^{(0)}  &=&  \frac{1}{{4\pi }}\int\limits_{\mathbb{R}^3 }{\frac{{u_2\,\mathcal{G}^0\left(x_1-u_1, x_2-u_2,x_3-u_3\right) }}
{{\left[u_1^2+u_2^2+u_3^2 \right]^{\frac{3}{2}} }}} \prod\limits_{j = 1}^3 {du_j }
\nonumber\\
&  = &\pi \left\{ {ac\sin \left(\pi x_2\right) \sin \left(\pi x_3\right)-bc\cos \left(\pi x_1\right) \cos \left(\pi x_2\right) } \right\}\qquad
\end{eqnarray}
%
\begin{eqnarray}
\label{A14}
g_3^{(0)}  &=&  \frac{1}{{4\pi }}\int\limits_{\mathbb{R}^3 }{\frac{{u_3\,\mathcal{G}^0\left(x_1-u_1, x_2-u_2,x_3-u_3\right) }}
{{\left[u_1^2+u_2^2+u_3^2 \right]^{\frac{3}{2}} }}} \prod\limits_{j = 1}^3 {du_j } 
\nonumber\\
&  = & \pi \left\{ {ab\sin\left( \pi x_1\right) \sin \left(\pi x_3\right)-ac\cos\left( \pi x_2\right) \cos \left(\pi x_3\right)} \right\}\qquad
\end{eqnarray}
resulting in $\mathcal{U}_i\left( x,t \right) \equiv 0$.  $v\left( {x,t} \right)$ is obtained from $\left(\ref{3.12}\right)$:
%
\begin{eqnarray}
\label{A15}
\frac{v_1}{v_r} = \left\{ {a\sin \left( {\pi x_3 } \right) - c\cos \left( {\pi x_2 } \right)} \right\}e^{ - \pi ^2 \kappa t} 
\end{eqnarray}
%
\begin{eqnarray}
\label{A16}
\frac{v_2}{v_r}   = \left\{b\sin \left( {\pi x_1 } \right)- {a\cos \left( {\pi x_3 } \right)} \right\}e^{ - \pi ^2 \kappa t} 
\end{eqnarray}
%
\begin{eqnarray}
\label{A17}
\frac{v_3}{v_r}   = \left\{ {c\sin \left( {\pi x_2 } \right) - b\cos \left( {\pi x_1 } \right)} \right\}e^{ - \pi ^2 \kappa t} 
\end{eqnarray}
A unit periodicity is observed for the solution according to the chosen initial conditions $\left(\ref{A8}\right)$-$\left(\ref{A10}\right)$; no velocity component is dependent on all three spatial coordinates. The integrations were carried out using the integral identities $\left(\ref{B1}\right)$-$\left(\ref{B3}\right)$ in Appendix B. Taking into account the unit periodicity of the velocity field $\left(m=1\right)$, integrating $\left(\ref{3.13}\right)$, we  obtain the pressure $p\left( {x,t} \right)$):
\begin{eqnarray}
\label{A18}
p \!= \! - \rho e^{ - 2\pi ^2 \kappa t}\! \left[ {bc\cos \left( {\pi x_1 } \right)\sin \left( {\pi x_2 } \right) + } {ab\cos \left( {\pi x_3 } \right)\sin \left( {\pi x_1 } \right) + ac\cos \left( {\pi x_2 } \right)\sin \left( {\pi x_3 } \right)} \right]\quad\quad
\end{eqnarray}
\section{Integral identities used in this paper}
%
\begin{eqnarray}
\label{B1}
\int\limits_{ - \infty }^\infty  {{e^{ - \frac{{{{\left( {x - u} \right)}^2}}}{{4\tau }}}}} du = 2\sqrt {\pi \tau }
\end{eqnarray}
%
\begin{eqnarray}
\label{C2}
 \int\limits_{ - \infty }^\infty  {\sin \left( {\alpha u} \right){e^{ - \frac{{{{\left( {x - u} \right)}^2}}}{{4\tau }}}}} du = 2\sqrt {\pi \tau } {e^{ - {\alpha ^2}\tau }}\sin \left( {\alpha x} \right) 
\end{eqnarray}
%
\begin{eqnarray}
\label{B3}
\int\limits_{ - \infty }^\infty  {\cos \left( {\alpha u} \right){e^{ - \frac{{{{\left( {x - u} \right)}^2}}}{{4\tau }}}}} du = 2\sqrt {\pi \tau } {e^{ - {\alpha ^2}\tau }}\cos \left( {\alpha x} \right)
\end{eqnarray}
%
\begin{eqnarray}
\label{B4}
\int\limits_{-\infty}^\infty\int\limits_{-\infty}^\infty\int\limits_{-\infty}^\infty {\frac{{u_1\cos\{\beta\left(x_1-u_1\right)\}}}{{\left[ u_1^2+u_2^2+u_3^2 \right]^{\frac{3}{2}} }}}\prod\limits_{j = 1}^3 {du_j }=\frac{4\pi}{\beta}\sin\left(\beta x_1\right)
\end{eqnarray}
%
\begin{eqnarray}
\label{B5}
\int\limits_{-\infty}^\infty\int\limits_{-\infty}^\infty\int\limits_{-\infty}^\infty \!{\frac{{u_1\cos\{\beta\left(x_1-u_1\right)\}\cos\{\beta\left(x_2-u_2\right)\}}}{{\left[ u_1^2+u_2^2+u_3^2 \right]^{\frac{3}{2}} }}}\prod\limits_{j = 1}^3 {du_j }=\frac{2\pi}{\beta}\sin\left(\beta x_1\right)\cos\left(\beta x_2\right)\quad
\end{eqnarray}
%
\begin{eqnarray}
\label{B6}
\int\limits_{-\infty}^\infty\!\int\limits_{-\infty}^\infty\!\int\limits_{-\infty}^\infty \!{\frac{{u_1\sin\{\beta\left(x_1-u_1\right)\}\cos\{\beta\left(x_2-u_2\right)\}}}{{\left[ u_1^2+u_2^2+u_3^2 \right]^{\frac{3}{2}} }}}\prod\limits_{j = 1}^3 {du_j }=-\frac{2\pi}{\beta}\cos\left(\beta x_1\right)\cos\left(\beta x_2\right)\quad
\end{eqnarray}
\begin{eqnarray}
\label{B7}
\int\limits_{-\infty}^\infty\!\int\limits_{-\infty}^\infty\!\int\limits_{-\infty}^\infty\! {\frac{{u_1\cos\{\beta\left(x_1-u_1\right)\}\sin\{\beta\left(x_2-u_2\right)\}}}{{\left[ u_1^2+u_2^2+u_3^2 \right]^{\frac{3}{2}} }}}\prod\limits_{j = 1}^3 {du_j }=\frac{2\pi}{\beta}\sin\left(\beta x_1\right)\sin\left(\beta x_2\right)\quad
\end{eqnarray}
%
\begin{eqnarray}
\label{B8}
\int\limits_{-\infty}^\infty\!\int\limits_{-\infty}^\infty\!\int\limits_{-\infty}^\infty\! {\frac{{u_1\sin\{\beta\left(x_1-u_1\right)\}\sin\{\beta\left(x_2-u_2\right)\}}}{{\left[ u_1^2+u_2^2+u_3^2 \right]^{\frac{3}{2}} }}}\prod\limits_{j = 1}^3 {du_j }=-\frac{2\pi}{\beta}\cos\left(\beta x_1\right)\sin\left(\beta x_2\right)\quad
\end{eqnarray}
%
\begin{eqnarray}
\label{B9}
\!\!\!\!\int\limits_{-\infty}^\infty\!\int\limits_{-\infty}^\infty\!\int\limits_{-\infty}^\infty \!\int\limits_{-\infty}^\infty\!{\frac{{u_1\sin\{\beta\left(x_1-u_1\right)\}\sin\{\beta\left(x_3-u_3\right)\}}}{{\left[ u_1^2+u_2^2+u_3^2+u_4^2 \right]^2 }}}\prod\limits_{j = 1}^4 {du_j }=-\frac{2\pi^2}{\beta}\cos\left(\beta x_1\right)\sin\left(\beta x_3\right)\quad
\end{eqnarray}
where $\alpha$ and $\beta$ are real constants.
\section{Derived expressions for $\mathcal{V}_1^0\left( {x_2,x_3,x_4; \xi_1, \xi_2, \xi_3, \xi_4} \right)$ and  $\mathcal{W}_1^0\left( {x_1,x_2,x_3,x_4; \xi_1, \xi_2, \xi_3, \xi_4} \right)$ for the four dimensional solution}
%
\protect\tiny{\begin{eqnarray}
\label{C1}
&&\frac{\mathcal{V}_1^0\left( {x_2,x_3,x_4; \xi_1, \xi_2, \xi_3, \xi_4} \right)}{v_r^2}=\nonumber\\
&-&\frac{\alpha }{2}\sin { \left(\xi_1-\xi_4 \right) }\sin { \left(\alpha x_2+\xi_1 \right) } \cos {\left(\alpha x_3+\xi_2\right) }\sin { \left(\alpha x_4+\xi_3\right) }  \sin{\left(\alpha x_2+\xi_2\right) }\sin { \left(\alpha x_3+\xi_3\right) } \cos { \left(\alpha x_4+\xi_4 \right) } -\nonumber\\
&-&\frac{\alpha }{2}\sin  {\left(\xi_2-\xi_4\right) } \sin { \left(\alpha x_2+\xi_1 \right) } \cos {\left(\alpha x_3+\xi_2\right) }\sin { \left(\alpha x_4+\xi_3\right) } \cos  { \left(\alpha x_3+\xi_4 \right) } \cos  {\left( \alpha x_2+\xi_3\right) } \cos  { \left(\alpha x_4+\xi_1 \right) }+\nonumber\\
&+&\frac{\alpha }{2}\sin  {\left( \xi_2-\xi_1\right) }\sin  {\left( \alpha x_2+\xi_2\right) } \cos  { \left(\alpha x_4+\xi_4 \right) } \sin  {\left( \alpha x_3+\xi_3\right) } \cos  { \left(\alpha x_3+\xi_4 \right) } \cos  {\left( \alpha x_2+\xi_3\right) } \cos  { \left(\alpha x_4+\xi_1 \right) }+\nonumber\\
&+&\frac{\alpha }{2}\cos { \left(\xi_1-\xi_3 \right) } \sin { \left(\alpha x_3+\xi_1 \right) } \cos {\left(\alpha x_4+\xi_2\right) } \cos { \left(\alpha x_2+\xi_4 \right) } \cos {\left(\alpha x_2+\xi_2\right) }\cos { \left(\alpha x_3+\xi_3\right) } \cos { \left(\alpha x_4+\xi_4 \right) } -\nonumber\\
&-&\frac{\alpha }{2}\sin { \left(\xi_1-\xi_4 \right) } \sin  {\left( \alpha x_3+\xi_2\right) }\sin  {\left( \alpha x_4+\xi_3\right) }\cos  {\left( \alpha x_2+\xi_1\right) }\cos {\left(\alpha x_2+\xi_2\right) }\cos { \left(\alpha x_3+\xi_3\right) } \cos { \left(\alpha x_4+\xi_4 \right) } +\nonumber\\
&+&\frac{\alpha }{2}\cos  {\left( \xi_2-\xi_3\right) } \sin { \left(\alpha x_3+\xi_1 \right) } \cos {\left(\alpha x_4+\xi_2\right) } \cos { \left(\alpha x_2+\xi_4 \right) } \sin  { \left(\alpha x_3+\xi_4 \right) } \sin  {\left( \alpha x_2+\xi_3\right) }\cos  {\left( \alpha x_4+\xi_1\right) } -\nonumber\\
&-&\frac{\alpha }{2}\sin  {\left( \xi_2-\xi_4\right) }\sin  {\left( \alpha x_3+\xi_2\right) }  \sin  {\left( \alpha x_4+\xi_3\right) }\cos  {\left( \alpha x_2+\xi_1\right) } \sin  { \left(\alpha x_3+\xi_4 \right) } \sin  {\left( \alpha x_2+\xi_3\right) }\cos  {\left( \alpha x_4+\xi_1\right) }-\nonumber\\
&-&\frac{\alpha }{2}\sin { \left(\xi_1-\xi_2 \right) } \sin { \left(\alpha x_4+\xi_1 \right) }\sin { \left(\alpha x_2+\xi_3\right) } \cos { \left(\alpha x_3+\xi_4 \right) } \cos {\left(\alpha x_2+\xi_2\right) }\sin { \left(\alpha x_3+\xi_3\right) } \sin { \left(\alpha x_4+\xi_4 \right) } +\nonumber\\
&+&\frac{\alpha }{2}\cos{ \left(\xi_1-\xi_3 \right) }  \sin  {\left( \alpha x_4+\xi_2\right) } \cos  { \left(\alpha x_2+\xi_4 \right) } \cos  {\left( \alpha x_3+\xi_1\right) }\cos {\left(\alpha x_2+\xi_2\right) }\sin { \left(\alpha x_3+\xi_3\right) } \sin { \left(\alpha x_4+\xi_4 \right) }-\nonumber\\
&-&\frac{\alpha }{2}\cos  {\left(\xi_2-\xi_3\right) } \sin  {\left( \alpha x_4+\xi_2\right) } \cos  { \left(\alpha x_2+\xi_4 \right) } \cos  {\left( \alpha x_3+\xi_1\right) } \cos  { \left(\alpha x_3+\xi_4 \right) } \sin  {\left( \alpha x_2+\xi_3\right) }\sin {\left( \alpha x_4+\xi_1\right) }
\end{eqnarray}}
\nopagebreak[4]
\begin{eqnarray}
\label{C2}
&&\frac{\mathcal{W}_1^0\left( {x_1,x_2,x_3,x_4; \xi_1, \xi_2, \xi_3, \xi_4} \right)}{{v_r^2}}=\frac{\alpha }{16}\sin \left[2{ \left(\alpha x_1+\xi_1 \right) }\right]+ \frac{\alpha }{16}\sin \left[2{ \left(\alpha x_1+\xi_1 \right) }\right]\cos\left[2 {\left(\alpha x_2+\xi_2 \right) }\right]-\nonumber\\
&-&\frac{\alpha }{16}\sin \left[2{ \left(\alpha x_1+\xi_1 \right) }\right]\cos\left[2{ \left(\alpha x_3+\xi_3\right) }\right]-\frac{\alpha }{16}\sin \left[2{ \left(\alpha x_1+\xi_1 \right) }\right]\cos\left[2 {\left(\alpha x_2+\xi_2 \right) }\right]\cos\left[2{ \left(\alpha x_3+\xi_3\right) }\right]+\nonumber\\
&+& \frac{\alpha }{16}\sin \left[2{ \left(\alpha x_1+\xi_1 \right) }\right]\cos\left[2{ \left(\alpha x_4+\xi_4\right) }\right]+ \frac{\alpha }{16}\sin \left[2{ \left(\alpha x_1+\xi_1 \right) }\right]\cos\left[2 {\left(\alpha x_2+\xi_2 \right) }\right]\cos\left[2{ \left(\alpha x_4+\xi_4\right) }\right]-\nonumber\\
&-& \frac{\alpha }{16}\sin \left[2{ \left(\alpha x_1+\xi_1 \right) }\right]\cos\left[2{ \left(\alpha x_3+\xi_3\right) }\right]\cos\left[2{ \left(\alpha x_4+\xi_4\right) }\right]- \nonumber\\
&-&\frac{\alpha }{16}\sin \left[2{ \left(\alpha x_1+\xi_1 \right) }\right]\cos\left[2 {\left(\alpha x_2+\xi_2 \right) }\right]\cos\left[2{ \left(\alpha x_3+\xi_3\right) }\right]\cos\left[2{ \left(\alpha x_4+\xi_4\right) }\right]+\nonumber\\
&+&\frac{\alpha }{16} \sin\left[ 2 {\left( \alpha x_1+\xi_2\right) }\right] -\frac{\alpha }{16} \sin\left[ 2 {\left( \alpha x_1+\xi_2\right) }\right]\cos\left[2{ \left(\alpha x_2+\xi_3\right) }\right]+\frac{\alpha }{16} \sin\left[ 2 {\left( \alpha x_1+\xi_2\right) }\right]\cos\left[2{ \left(\alpha x_3+\xi_4\right) }\right]-\nonumber\\
&-&\frac{\alpha }{16} \sin\left[ 2 {\left( \alpha x_1+\xi_2\right) }\right]\cos\left[2{ \left(\alpha x_2+\xi_3\right) }\right]\cos\left[2{ \left(\alpha x_3+\xi_4\right) }\right]+\frac{\alpha }{16} \sin\left[ 2 {\left( \alpha x_1+\xi_2\right) }\right]\cos\left[2{ \left(\alpha x_4+\xi_1\right) }\right]-\nonumber\\
&-&\frac{\alpha }{16} \sin\left[ 2 {\left( \alpha x_1+\xi_2\right) }\right]\cos\left[2{ \left(\alpha x_2+\xi_3\right) }\right]\cos\left[2{ \left(\alpha x_4+\xi_1\right) }\right]+\frac{\alpha }{16} \sin\left[ 2 {\left( \alpha x_1+\xi_2\right) }\right]\cos\left[2{ \left(\alpha x_3+\xi_4\right) }\right]\cos\left[2{ \left(\alpha x_4+\xi_1\right) }\right]-\nonumber\\
&-&\frac{\alpha }{16} \sin\left[ 2 {\left( \alpha x_1+\xi_2\right) }\right]\cos\left[2{ \left(\alpha x_2+\xi_3\right) }\right]\cos\left[2{ \left(\alpha x_3+\xi_4\right) }\right]\cos\left[2{ \left(\alpha x_4+\xi_1\right) }\right]- \nonumber\\
&-&\alpha\sin  {\left( 2\alpha x_1+\xi_1+\xi_2\right) }\cos  { \left(\alpha x_3+\xi_4 \right) } \sin  {\left( \alpha x_2+\xi_3\right) }\cos  {\left( \alpha x_4+\xi_1\right) } \cos {\left(\alpha x_2+\xi_2\right) }\sin { \left(\alpha x_3+\xi_3\right) } \cos { \left(\alpha x_4+\xi_4 \right) }-\nonumber\\
&-&\frac{\alpha }{2}\sin { \left(2\alpha x_1+\xi_1+\xi_4 \right) } \sin { \left(\alpha x_2+\xi_1 \right) } \cos {\left(\alpha x_3+\xi_2\right) }\sin { \left(\alpha x_4+\xi_3\right) } \sin{\left(\alpha x_2+\xi_2\right) }\sin { \left(\alpha x_3+\xi_3\right) } \cos { \left(\alpha x_4+\xi_4 \right) } -\nonumber\\
&-&\frac{\alpha }{2}\sin  {\left( 2\alpha x_1+\xi_2+\xi_4 \right) }\sin { \left(\alpha x_2+\xi_1 \right) } \cos {\left(\alpha x_3+\xi_2\right) }\sin { \left(\alpha x_4+\xi_3\right) }  \cos  { \left(\alpha x_3+\xi_4 \right) } \cos  {\left( \alpha x_2+\xi_3\right) } \cos  { \left(\alpha x_4+\xi_1 \right) }+\nonumber\\
&+&\frac{\alpha }{2}\sin  {\left(2 \alpha x_1+\xi_1+\xi_2\right) }\sin  {\left( \alpha x_2+\xi_2\right) } \cos  { \left(\alpha x_4+\xi_4 \right) } \sin  {\left( \alpha x_3+\xi_3\right) } \cos  { \left(\alpha x_3+\xi_4 \right) } \cos  {\left( \alpha x_2+\xi_3\right) } \cos  { \left(\alpha x_4+\xi_1 \right) }+\nonumber\\
&+&\frac{\alpha }{16}\sin \left[2{ \left(\alpha x_1+\xi_1 \right) }\right]-\frac{\alpha }{16}\sin \left[2{ \left(\alpha x_1+\xi_1 \right) }\right]\cos\left[2 {\left(\alpha x_2+\xi_2 \right) }\right]- \frac{\alpha }{16}\sin \left[2{ \left(\alpha x_1+\xi_1 \right) }\right]\cos\left[2{ \left(\alpha x_3+\xi_3\right) }\right]+\nonumber\\
&+& \frac{\alpha }{16}\sin \left[2{ \left(\alpha x_1+\xi_1 \right) }\right]\cos\left[2 {\left(\alpha x_2+\xi_2 \right) }\right]\cos\left[2{ \left(\alpha x_3+\xi_3\right) }\right]+ \frac{\alpha }{16}\sin \left[2{ \left(\alpha x_1+\xi_1 \right) }\right]\cos\left[2{ \left(\alpha x_4+\xi_4\right) }\right]-\nonumber\\
&-& \frac{\alpha }{16}\sin \left[2{ \left(\alpha x_1+\xi_1 \right) }\right]\cos\left[2 {\left(\alpha x_2+\xi_2 \right) }\right]\cos\left[2{ \left(\alpha x_4+\xi_4\right) }\right]- \frac{\alpha }{16}\sin \left[2{ \left(\alpha x_1+\xi_1 \right) }\right]\cos\left[2{ \left(\alpha x_3+\xi_3\right) }\right]\cos\left[2{ \left(\alpha x_4+\xi_4\right) }\right]+\nonumber\\
&+&\frac{\alpha }{16}\sin \left[2{ \left(\alpha x_1+\xi_1 \right) }\right]\cos\left[2 {\left(\alpha x_2+\xi_2 \right) }\right]\cos\left[2{ \left(\alpha x_3+\xi_3\right) }\right]\cos\left[2{ \left(\alpha x_4+\xi_4\right) }\right]-\nonumber\\
&-&\frac{\alpha }{2}\cos { \left(2\alpha x_1+\xi_1+\xi_3\ \right) } \sin { \left(\alpha x_3+\xi_1 \right) } \cos {\left(\alpha x_4+\xi_2\right) } \cos { \left(\alpha x_2+\xi_4 \right) }\cos {\left(\alpha x_2+\xi_2\right) }\cos { \left(\alpha x_3+\xi_3\right) } \cos { \left(\alpha x_4+\xi_4 \right) } -\nonumber\\
&-&\frac{\alpha }{2}\sin { \left(2\alpha x_1+\xi_1+\xi_4 \right) } \sin  {\left( \alpha x_3+\xi_2\right) } \sin  {\left( \alpha x_4+\xi_3\right) }\cos  {\left( \alpha x_2+\xi_1\right) }  \cos {\left(\alpha x_2+\xi_2\right) }\cos { \left(\alpha x_3+\xi_3\right) } \cos { \left(\alpha x_4+\xi_4 \right) } -\nonumber\\
&-&\frac{\alpha }{2}\cos  {\left( 2\alpha x_1+\xi_2+\xi_3\right) }\sin { \left(\alpha x_3+\xi_1 \right) } \cos {\left(\alpha x_4+\xi_2\right) } \cos { \left(\alpha x_2+\xi_4 \right) }\sin  { \left(\alpha x_3+\xi_4 \right) } \sin  {\left( \alpha x_2+\xi_3\right) }\cos  {\left( \alpha x_4+\xi_1\right) } -\nonumber\\
&-&\frac{\alpha }{2}\sin  {\left( 2\alpha x_1+\xi_2+\xi_4\right) }\sin  {\left( \alpha x_3+\xi_2\right) }  \sin  {\left( \alpha x_4+\xi_3\right) }\cos  {\left( \alpha x_2+\xi_1\right) }\sin  { \left(\alpha x_3+\xi_4 \right) } \sin  {\left( \alpha x_2+\xi_3\right) }\cos  {\left( \alpha x_4+\xi_1\right) }+\nonumber\\
&+&\frac{\alpha }{16}\sin \left[2{ \left(\alpha x_1+\xi_2 \right) }\right]-\frac{\alpha }{16}\sin \left[2{ \left(\alpha x_1+\xi_2 \right) }\right]\cos\left[2 {\left(\alpha x_4+\xi_1 \right) }\right]- \frac{\alpha }{16}\sin \left[2{ \left(\alpha x_1+\xi_2 \right) }\right]\cos\left[2{ \left(\alpha x_2+\xi_3\right) }\right]+\nonumber\\
&+& \frac{\alpha }{16}\sin \left[2{ \left(\alpha x_1+\xi_2 \right) }\right]\cos\left[2 {\left(\alpha x_2+\xi_3 \right) }\right]\cos\left[2{ \left(\alpha x_4+\xi_1\right) }\right]+ \frac{\alpha }{16}\sin \left[2{ \left(\alpha x_1+\xi_2 \right) }\right]\cos\left[2{ \left(\alpha x_3+\xi_4\right) }\right]-\nonumber\\
&-& \frac{\alpha }{16}\sin \left[2{ \left(\alpha x_1+\xi_2 \right) }\right]\cos\left[2 {\left(\alpha x_3+\xi_4 \right) }\right]\cos\left[2{ \left(\alpha x_4+\xi_1\right) }\right]- \frac{\alpha }{16}\sin \left[2{ \left(\alpha x_1+\xi_2 \right) }\right]\cos\left[2{ \left(\alpha x_2+\xi_3\right) }\right]\cos\left[2{ \left(\alpha x_3+\xi_4\right) }\right]+ \nonumber\\
&+&\frac{\alpha }{16}\sin \left[2{ \left(\alpha x_1+\xi_2 \right) }\right]\cos\left[2 {\left(\alpha x_2+\xi_3 \right) }\right]\cos\left[2{ \left(\alpha x_3+\xi_4\right) }\right]\cos\left[2{ \left(\alpha x_4+\xi_1\right) }\right]-\nonumber\\
&-&\frac{\alpha }{2} \sin { \left(2\alpha x_1+\xi_1+\xi_2 \right) }\sin { \left(\alpha x_4+\xi_1 \right) } \sin { \left(\alpha x_2+\xi_3\right) } \cos { \left(\alpha x_3+\xi_4 \right) } \cos {\left(\alpha x_2+\xi_2\right) }\sin { \left(\alpha x_3+\xi_3\right) } \sin { \left(\alpha x_4+\xi_4 \right) } +\nonumber\\
&-&\frac{\alpha }{2}\cos { \left(\alpha x_1+\xi_1+\xi_3 \right) } \sin  {\left( \alpha x_4+\xi_2\right) } \cos  { \left(\alpha x_2+\xi_4 \right) }\cos  {\left( \alpha x_3+\xi_1\right) }\cos {\left(\alpha x_2+\xi_2\right) }\sin { \left(\alpha x_3+\xi_3\right) } \sin { \left(\alpha x_4+\xi_4 \right) }+\nonumber\\
&+&\frac{\alpha }{2}\cos {\left( 2\alpha x_1+\xi_2+\xi_3\right) } \sin  {\left( \alpha x_4+\xi_2\right) } \cos  { \left(\alpha x_2+\xi_4 \right) } \cos  {\left( \alpha x_3+\xi_1\right) }\cos  { \left(\alpha x_3+\xi_4 \right) } \sin  {\left( \alpha x_2+\xi_3\right) }\sin {\left( \alpha x_4+\xi_1\right) }
\end{eqnarray}
\bibliographystyle{chicago}
\bibliography{NS}
\end{document}